\documentclass[12pt,reqno]{amsart}
\usepackage{amsmath,amsthm,amssymb,amsfonts,amscd}
\usepackage{mathrsfs}
\usepackage{bbm}
\usepackage{bbding}

\usepackage{bm}
\usepackage[backref=page]{hyperref}

\usepackage{geometry}
\usepackage{color}
\usepackage{xcolor}

\usepackage{picture,epic}
\usepackage{tikz}

\numberwithin{equation}{section}

\theoremstyle{plain}
\newtheorem{theorem}{Theorem}[section]
\newtheorem{lemma}[theorem]{Lemma}

\theoremstyle{definition}

\theoremstyle{remark}
\newtheorem{remark}[theorem]{Remark}

\renewcommand{\Re}{\operatorname{Re}}
\renewcommand{\Im}{\operatorname{Im}}

\newcommand{\sym}{\operatorname{sym}}

\newcommand{\GL}{\operatorname{GL}}
\newcommand{\SL}{\operatorname{SL}}

\newcommand{\dd}{\mathrm{d}}

\makeatletter
\def\@tocline#1#2#3#4#5#6#7{\relax
  \ifnum #1>\c@tocdepth 
  \else
    \par \addpenalty\@secpenalty\addvspace{#2}%
    \begingroup \hyphenpenalty\@M
    \@ifempty{#4}{%
      \@tempdima\csname r@tocindent\number#1\endcsname\relax
    }{%
      \@tempdima#4\relax
    }%
    \parindent\z@ \leftskip#3\relax \advance\leftskip\@tempdima\relax
    \rightskip\@pnumwidth plus4em \parfillskip-\@pnumwidth
    #5\leavevmode\hskip-\@tempdima
      \ifcase #1
       \or\or \hskip 1em \or \hskip 2em \else \hskip 3em \fi%
      #6\nobreak\relax
    \hfill\hbox to\@pnumwidth{\@tocpagenum{#7}}\par
    \nobreak
    \endgroup
  \fi}
\makeatother

\begin{document}

\title[Determination of $\GL(3)$ cusp forms]{Determination of $\GL(3)$ cusp forms by central values of quadratic twisted  $L$-functions}

\author{Shenghao Hua and Bingrong Huang}

\address{Data Science Institute and School of Mathematics \\ Shandong University \\ Jinan \\ Shandong 250100 \\China}
\email{huashenghao@vip.qq.com}
\email{brhuang@sdu.edu.cn}

\date{\today}

\begin{abstract}
  Let  $\phi$ and $\phi'$ be two  $\GL(3)$ Hecke--Maass cusp forms. In this paper, we prove  that $\phi=\phi'\textrm{ or }\widetilde{\phi'}$ if  there exists a nonzero constant $\kappa$ such that   $L(\frac{1}{2},\phi\otimes \chi_{8d})=\kappa L(\frac{1}{2},\phi'\otimes \chi_{8d})$ for all positive odd square-free positive  $d$. Here $\widetilde{\phi'}$ is dual form of $\phi'$ and $\chi_{8d}$ is the quadratic character $(\frac{8d}{\cdot})$.
  To prove this, we obtain asymptotic formulas for twisted first moment of central values of quadratic twisted  $L$-functions on $\GL(3)$, which will have many other applications.
\end{abstract}



\thanks{This work was supported by  the National Key R\&D Program of China (No. 2021YFA1000700),  NSFC (Nos. 12001314 and 12031008), and the Young Taishan Scholars Program.}

\maketitle

\section{Introduction} \label{sec:Intr}

Determining automorphic forms from central values of the twisted $L$-functions is a topic of much interest (see e.g. \cite{LR97,Luo99,CD05,MS15, KL20}).
It was first considered by Luo and Ramakrishnan \cite{LR97} for modular forms.
They showed that if two cuspidal normalized newforms $f$ and $g$ of weight $2k$  (resp. $2k'$) and level $N$ (resp. $N'$) have the property that
\begin{equation}\label{need}
L (\frac{1}{2}, f\otimes\chi_d ) = L (\frac{1}{2}, g\otimes\chi_d )
\end{equation}
for all quadratic characters $\chi_d$, then $k=k'$, $N=N'$ and $f=g$.
Chinta and Diaconu \cite{CD05} proved that self-dual $\GL(3)$ Hecke--Maass forms are determined by their quadratic twisted central $L$-values.  
They used the method of double Dirichlet series for the averaging process.
Recently, Kuan and Lesesvre \cite{KL20} generalized the analogous result to automorphic representations of $\GL(3,F)$ over number field which are self-contragredient.\footnote{They assume ``Hypothesis 1'' which is satisfied by self-contragredient forms.}
In this paper, we use a new method to give a general result on $\GL(3)$, without the assumptions of \cite{KL20}. Instead of using double Dirichlet series as in \cite{CD05} and \cite{KL20}, we introduce a twisted average of the central $L$-values and obtain its asymptotics for which we use a method based on Soundararajan's work \cite{SD00}.

Let $\phi$ be a Hecke--Maass cusp form of type $\nu=(\nu_1,\nu_2)$ for $\textrm{SL}(3,\mathbb{Z})$ with the normalized Fourier coefficients
 $A(m,n)$.
We have the conjugation relation $A(m,n)=\overline{A(n,m)}$, see \cite[Theorem 9.3.11]{GB06}.
Any real primitive character to the modulus $q$ must be of the form $\chi(n)=(\frac{d}{n})$ where $d$ is a fundamental discriminant \cite[Theorem 9.13]{MV07}, i.e., a product of pairwise coprime integers of the form $-4$, $\pm 8$, $(-1)^{\frac{p-1}{2}}p$ where $p$ is an odd prime.
There are two primitive characters to the modulus $q$ if $8\parallel q$ and only one otherwise. From \cite[Theorem 7.1.3]{GB06} we know
\begin{equation*}
\prod_{i=1}^3\Gamma(\frac{\frac{1}{2}-\gamma_i}{2})L(\frac{1}{2},\phi\otimes \chi_{8d})=\prod_{i=1}^3\Gamma(\frac{\frac{1}{2}
+\gamma_i}{2})L(\frac{1}{2},\widetilde{\phi}
\otimes \chi_{8d})
\end{equation*}
when $d$ is positive where $\widetilde{\phi}$ is the dual form of $\phi$, $\gamma_1=1-2\nu_1-\nu_2$, $\gamma_2=\nu_1-\nu_2$, and $\gamma_3=-1+\nu_1+2\nu_2$ are the Langlands parameters of $\phi$.
By unitarity and the standard Jacquet--Shalika bounds, the Langlands parameters of an arbitrary irreducible representation $\pi\subseteq L^2(\textrm{SL}(3,\mathbb{Z})\backslash \mathbb{H}^3)$ must satisfy $\sum_{i=1}^3\gamma_i=0$ and $\{-\gamma_i\}_{i=1}^3=\{\overline{\gamma_i}\}_{i=1}^3$.
Let
\begin{equation*}
  \kappa_{\phi,\phi'} = \left\{ \begin{array}{ll}
                     1, & \textrm{if $\phi=\phi'$}, \\
                     \prod_{i=1}^3\frac{\Gamma(\frac{\frac{1}{2}+\gamma_i}{2})} {\Gamma(\frac{\frac{1}{2}-\gamma_i}{2})}, &  \textrm{if $\phi=\widetilde{\phi'}$}.
                   \end{array}\right.
\end{equation*}
Then we have
\begin{equation*}
L(\frac{1}{2},\phi\otimes \chi_{8d})= \kappa_{\phi,\phi'} L(\frac{1}{2},\phi'\otimes \chi_{8d}),
\end{equation*}
for all positive fundamental discriminants $8d$.
But we don't know if the converse conclusion is true, i.e. if for normalized $\phi$ and $\phi'$,  there exists a nonzero constant $\kappa$ such that $L(\frac{1}{2},\phi\otimes \chi_{8d})=\kappa L(\frac{1}{2},\phi'\otimes \chi_{8d})$ for all positive $d$,
is it then true that we have$\phi=\phi'\textrm{ or }\widetilde{\phi'}$?
Our main result in this paper is as follows.

\begin{theorem}\label{thm1}
Let $\phi$ and $\phi'$ be two normalized Hecke--Maass cusp forms of $\operatorname{SL}(3,\mathbb{Z})$.
Fix an integer $M$ coprime to $3,5,7,11$.
If there exists a nonzero constant $\kappa$ such that
\begin{equation}\label{eq1}
L(\frac{1}{2},\phi\otimes \chi_{8d})=\kappa L(\frac{1}{2},\phi'\otimes \chi_{8d})
\end{equation}
hold for all positive odd square-free integers $d$ coprime to $M$, then we have $\phi=\phi' \textrm{ or }\widetilde{\phi'}$.
If we further assume $\prod_{i=1}^3\Gamma(\frac{\frac{1}{2}-\gamma_i}{2})\notin \mathbb{R}$, then we have $\phi=\phi'$ if and only if  $\kappa=1$, and $\phi=\widetilde{\phi'}$ if and only if $\kappa=\prod_{i=1}^3  \frac{\Gamma(\frac{\frac{1}{2}+\gamma_i}{2})} { \Gamma(\frac{\frac{1}{2}-\gamma_i}{2})}$.
\end{theorem}

\begin{remark}\label{rm1}
Let $\{\gamma_i\}$ be the Langlands parameters of $\phi$.
If $\prod_{i=1}^3\Gamma(\frac{\frac{1}{2}-\gamma_i}{2})\notin \mathbb{R}$, then $\phi$ is not self-dual. Our method also works for $\GL(3)$ forms of any fixed level.
\end{remark}


We prove Theorem \ref{thm1} by using the following Theorem \ref{thm2} on twisted first moment of central values of quadratic twisted $\GL(3)$ $L$-functions.
Our arguments combine ideas from Soundararajan \cite{SD00} and Chinta--Diaconu \cite{CD05}.
To extract the relevant information from the main term in Theorem \ref{thm2}, we will  use Lemma \ref{lm5.1} below.
This argument is different from Chinta--Diaconu's, where they used the fact that certain rational function is monotone for real variable (their Lemma 5.1), which may be special to the self-contragredient case since $A(p,1)$ are complex in the non self-contragredient case.

For an automorphic representation $\pi$ of $\GL(3,\mathbb{A_Q})$, $L(s,\sym^2\pi)$ has a simple pole at $s=1$ if and only if $\pi$ is the Gelbart--Jacquet lift \cite{GJ76} of an automorphic representation on $\GL(2,\mathbb{A_Q})$ with trivial central character \cite{GRS99}, i.e., it is a self-contragredient cuspidal automorphic representation \cite{RM14}.
So $\phi$ is self-dual if and only if $L(s,\sym^2 \phi)$ has a simple pole at $s=1$, which is equivalent to $L(s,\sym^2 \widetilde{\phi})$ has a simple pole at $s=1$.

We denote $\theta_3$ be the the least common upper bound of power of $p$ for $|A(p,1)|$, i.e., $|A(p,1)|\leq 3 p^{\theta_3}$ for all prime $p$.
The Generalized Ramanujan Conjecture implies that $\theta_3=0$, and from Kim--Sarnak \cite[Appendix 2]{Kim2003} we know $\theta_3\leq\frac{5}{14}$.

Let $\Phi$ be any smooth nonnegative Schwarz class function supported in the interval $(1,2)$.
For any integer $\nu \geq 0$ we define
\begin{equation*}
\Phi_{(\nu)}=\max_{0\leq j \leq \nu}\int_1^2|\Phi^{(j)}(t)| \dd t.
\end{equation*}
For any complex number $w$, we define
\begin{equation*}
\check{\Phi}(w)=\int_0^\infty \Phi(y)y^w \dd y,
\end{equation*}
so $\check{\Phi}(w)$ is holomorphic.
Integrating by parts $\nu$ times, we have
\begin{equation*}
\check{\Phi}(w)=\frac{1}{(w+1)\dots(w+\nu)}\int_0^\infty\Phi^{(\nu)}(y)y^{w+\nu}\dd y,
\end{equation*}
thus for $\Re (w)>-1$ we have
\begin{equation*}
|\check{\Phi}(w)|\ll_\nu \frac{2^{\textrm{Re}(w)}}{|w+1|^{\nu}}\Phi_{(\nu)}.
\end{equation*}

\begin{theorem}\label{thm2}
Let $\phi$ be a Hecke--Maass cusp form for $\operatorname{SL}(3,\mathbb{Z})$ with normalized Fourier coefficients $A(m,n)$.
For sufficiently large $X>0$, arbitrarily small $\varepsilon >0$, and any odd integer $l\ll X^{\frac{1}{15}-\varepsilon}$, if $\phi$ is not self-dual then we have
\begin{multline*}
\sum_{2\nmid d}\nolimits^\flat \chi_{8d}(l) L(\frac{1}{2},\phi\otimes \chi_{8d})\Phi(\frac{d}{X})
=\frac{2\check{\Phi}(0)X}{3\zeta(2)\sqrt{l_1}}
\prod_{p\mid l}\frac{p}{p+1}
 \bigg(G_\phi(l)L^{\{2\}}(1,\sym^2 \phi) \\ +\prod_{i=1}^3\frac{\Gamma(\frac{\frac{1}{2}
 +\gamma_i}{2})}{\Gamma(\frac{\frac{1}{2}
 -\gamma_i}{2})}
 \bar{G}_\phi(l)L^{\{2\}}(1,\sym^2 \widetilde{\phi})\bigg)
+O_\Phi ( l^{\frac{3}{4}+\varepsilon}X^{\frac{19}{20}+\varepsilon});
\end{multline*}
and if $\phi$ is self-dual then we have
\begin{multline*}
\sum_{2\nmid d}\nolimits^\flat \chi_{8d}(l) L(\frac{1}{2},\phi\otimes \chi_{8d})\Phi(\frac{d}{X})
= \frac{\lim_{s\to 1} (s-1)L^{\{2\}}(s,\sym^2 \phi)\check{\Phi}(0)}{\zeta(2)\sqrt{l_1}}
\prod_{p\mid l}\frac{p}{p+1}X\\
\times \Big(G_\phi(l)
\log \frac{X}{l_1 ^{\frac{2}{3}}}+C_\phi(l)\Big)
+O_\Phi ( l^{\frac{3}{4}+\varepsilon}
X^{\frac{19}{20}+\varepsilon}),
\end{multline*}
where $l=l_1l_2^2$ with $l_1$ is square-free, $\sum_{2\nmid d}\nolimits^\flat$ means summing  over positive odd square-free $d$,
$L^{\{2\}}(s, \sym ^2 \phi )=L(s, \sym ^2 \phi)/L_{2}(s, \sym ^2 \phi)$
and $L_{2}(s, \sym ^2 \phi)$ is the local $L$-function in $2$-place,
$C_\phi(l)$ is defined as in
\eqref{Cphi}, 
$G_\phi(l)=\prod_{\textrm{odd prime }p}G_{\phi,p}(l)$ with
\begin{equation}
G_{\phi,p}(l)=\left\{\begin{aligned}
&(A(p,1)+\frac{1}{p})
(1-\frac{A(1,p)}{p}
+\frac{A(p,1)}{p^{2}}-\frac{1}{p^{3}}),
&\textrm{if } p\mid l_1,\\
&(1+\frac{A(1,p)}{p})(1-\frac{A(1,p)}{p}
     +\frac{A(p,1)}{p^{2}}-\frac{1}{p^{3}}),
     &\textrm{if } p\mid l_2,p\nmid l_1,\\
&(1-\frac{p^{2}A(p,1)^2-
p^{2}A(1,p)-pA(1,p)^2+2pA(p,1)+1}
{(p+1)p^{2}(p+A(1,p))})
\\
&\hskip 30pt \times
(1+\frac{A(1,p)}{p})(1-\frac{A(1,p)}{p}
     +\frac{A(p,1)}{p^{2}}-\frac{1}{p^{3}}),
&\textrm{if } p\nmid l.
\end{aligned}
\right.
\end{equation}

Moreover,
there exists $c_\phi =3^a\times 5^b \times 7^2 \times 11^2$ with some $a,b\in\{1,2\}$
such that:
For $l=c_\phi l'$ with $(2c_\phi,l')=1$,
$G_\phi(l)= 0$ if and only if there exists a prime $13\leq p\mid l_1$ such that
$A(p,1)=-\frac{1}{p}$.
\end{theorem}

\begin{remark}
  We did not try to optimize the exponent $19/20$ as it is enough for our main theorem. One may compare our result with the cubic moment of quadratic Dirichlet $L$-functions (see e.g. \cite{SD00,Young,DGH,DW}). Our case is more complicated as $\phi$ is undecomposable.
  We will use Soundararajan's approach to prove Theorem \ref{thm2} which is based on the approximate functional equation and Poisson summation formula (Lemma \ref{lm:psf}).
  At present, one still can not  unconditionally prove an asymptotic formula (with power saving) for the fourth moment of quadratic Dirichlet $L$-functions. To extend our result to $\GL_n$ for $n > 3$ seems hard.
\end{remark}

\begin{remark}
In fact, our method of the proof also works for  the functions $n\mapsto d_3(n)=(1\star 1\star 1)(n)$ and $n\mapsto (\lambda_f \star 1)(n)$, which are the Hecke eigenvalues of certain non-cuspidal automorphic representations of $\GL(3,\mathbb{A}_\mathbb{Q})$, namely the isobaric representations $1\boxplus 1\boxplus 1$ and $\pi_f\boxplus 1$.
Here $\pi_f$ is a cuspidal automorphic representation of $\GL(2,\mathbb{A}_\mathbb{Q})$.
The method of the present paper can be generalized straightforwardly to show above results for an arbitrary irreducible automorphic representation $\pi \subseteq L^2(\SL(3,\mathbb{Z})\backslash \mathbb{H}^3)$.
\end{remark}

\begin{remark}
By using of the large sieve estimates for quadratic characters in \cite{HB95}, one may prove a nonvanishing result for central values of quadratic twisted $\GL(3)$ $L$-functions, i.e. there exist at least $O(X^{1/2-\varepsilon})$ fundamental discriminants $X\leq d\leq 2X$ for arbitrarily small $\varepsilon>0$ such that $L(\frac{1}{2},\phi\otimes \chi_{8d})\neq 0$.
\end{remark}

\begin{remark}
In \cite{HH22}, we give another application of Theorem \ref{thm2},
where we prove the conjectured order lower bounds for the $k$-th moments of central values of quadratic twisted self-dual $\GL(3)$ $L$-functions for all $k\geq 1$.
\end{remark}

To prove Theorem \ref{thm2}, we need the following Generalized Ramanujan Conjecture on average for a special sequence of Fourier coefficients of $\phi$, which is closely related to the symmetric square lift of $\phi$.
This may have its own interest.

\begin{theorem}\label{thm3}
Let $\phi$ be a fixed normalized Hecke--Maass cusp form for $\operatorname{SL}(3,\mathbb{Z})$, and $A(m,n)$ be its Fourier coefficients.
For any $\varepsilon>0$ we have
\begin{equation*}
\sum_{n\leq X}|A(n^2,1)|\ll_{\phi,\varepsilon} X^{1+\varepsilon}.\label{thm3eq1}
\end{equation*}
\end{theorem}

The rest of this paper is organized as follows.
In \S \ref{sec:preliminaries}, we introduce some notation and present some lemmas that we will need later.
In \S \ref{sec:first_moment}, we extend Soundararajan's method to prove Theorem \ref{thm2}.
In \S \ref{sec:det_form}, we prove Theorem \ref{thm1} by using Theorem \ref{thm2}.
Finally, in \S \ref{sec:FC}, we prove Theorem \ref{thm3} by using the Rankin--Selberg bounds on the averages of Fourier coefficients.

\section{Notation and preliminary results}\label{sec:preliminaries}
For any complex numbers sequence $\{f_n\}_{n=1}^\infty$ and smooth nonnegative Schwarz class function $\Phi$ supported in the interval $(1,2)$.
We define
\begin{equation*}
S(f_d;\Phi)=S_X(f_d;\Phi)=\frac{1}{X}\sum_{2\nmid d}\nolimits^\flat f_d\Phi(\frac{d}{X})=\frac{1}{X}\sum_{d \textrm{ odd}}\mu^2(d)f_d\Phi(\frac{d}{X}).
\end{equation*}

For real parameter $Y>1$ and we have $\mu^2(d)=M_Y(d)+R_Y(d)$ where
\begin{equation*}
M_Y(d)=\sum_{\substack{l^2\mid d \\ l \leq Y}}\mu(l),\textrm{ and }R_Y(d)=\sum_{\substack{l^2\mid d \\l > Y}}\mu(l).
\end{equation*}
Define
\begin{equation*}
S_M(f_d;\Phi)=S_{M,X,Y}(f_d;\Phi)=\frac{1}{X}\sum_{d \textrm{ odd}}M_Y(d)f_d\Phi(\frac{d}{X}),
\end{equation*}
and
\begin{equation*}
S_R(f_d;\Phi)=S_{R,X,Y}(f_d;\Phi)=\frac{1}{X}\sum_{d \textrm{ odd}}R_Y(d)f_d\Phi(\frac{d}{X}),
\end{equation*}
so $S(f_d;\Phi)=S_M(f_d;\Phi)+S_R(f_d;\Phi).$

Let $H(s)$ be any function which is holomorphic and bounded in the strip $-4<\Re (u)<4$, even, and normalized by $H(0)=1$.
From Kim--Sarnak \cite[Appendix 2]{Kim2003} we know
$$\max\{\Re(\gamma_i)\}\leq \frac{5}{14},$$
and from \cite{IK04} we know $L(s,\phi\otimes \chi_{8d})$ is entire and we have the following lemma.
\begin{lemma}[Approximate functional equation]  \label{lemma:AFE}
Let $d$ be a positive odd square-free integer.
Then we have
\begin{equation*}
  L(\frac{1}{2},\phi\otimes \chi_{8d})
  =\sum_{n=1}^\infty\bigg(A(n,1)V(n(\frac{\pi}{8d})^{\frac{3}{2}})
  +\prod_{i=1}^3\frac{\Gamma(\frac{
  \frac{1}{2}+\gamma_i}{2})} {\Gamma(\frac{\frac{1}{2}-\gamma_i}{2})}
  \overline{A(n,1)V(n(\frac{\pi}{ 8d})^{\frac{3}{2}})}\bigg)\frac{\chi_{8d}(n)}{\sqrt{n}},
\end{equation*}
where $V(y)$ is defined by
\begin{equation} \label{eqn:V}
V(y)=\frac{1}{2\pi i}\int_{(u)}\prod_{i=1}^3\frac{
\Gamma(\frac{s+\frac{1}{2}-\gamma_i}{2})} {\Gamma(\frac{\frac{1}{2}-\gamma_i}{2})}
y^{-s}H(s)\frac{\dd s}{s}
\end{equation}
with $u>0$.  Here, and in sequel, $\int_{(u)}$ stands for $\int_{u-i \infty}^{u+i \infty}$.
We can choose $H(s)=e^{s^2}$.
\end{lemma}

\begin{proof}
  See \cite[Theorem 5.3]{IK04}.
\end{proof}

\begin{lemma}\label{lm2.2}
The function $V$ is smooth on $[0,\infty)$. Let $h\in\mathbb{Z}_{\geq0}$.
For $y$ near $0$ and $v<\frac{1}{2}-\theta_3$ we have
\begin{equation*}
  y^h V^{(h)}(y)=\delta_h+O_h(y^{v}),
\end{equation*}
and for large $y$, any $A>0$, and any integer $h$,
\begin{equation*}
  y^h V^{(h)}(y)\ll_{h,A}  y^{-A}.
\end{equation*}
Here $\delta_0=1$, and $\delta_h=0$ if $h\geq1$.
\end{lemma}

\begin{proof}
  See \cite[Proposition 5.4]{IK04}.
\end{proof}

%
%

Let $n$ be an odd integer.
We define for all integers $k$
\begin{equation*}
G_k(n)=(\frac{1-i}{2}+(\frac{-1}{n})\frac{1+i}{2})\sum_{a\pmod{n}}(\frac{a}{n})e(\frac{ak}{n}),
\end{equation*}
and put
\begin{equation*}
\tau_k(n)=\sum_{a\pmod{n}}(\frac{a}{n})e(\frac{ak}{n})=(\frac{1+i}{2}+(\frac{-1}{n})\frac{1-i}{2})G_k(n).
\end{equation*}
Here $e(x)=\exp(2\pi ix)$.
If $n$ is square-free then $(\frac{\cdot}{n})$ is a primitive character with conductor $n$.
Here it is easy to see that $G_k(n)=(\frac{k}{n})\sqrt{n}$.
For our later work, we require knowledge of $G_k(n)$ for all odd $n$.

For fundamental discriminant $d$ we know Gauss sum of $\chi_d$ is $\tau(\chi_d)=\sqrt{d}$ where the square root is taken as its principal branch.

\begin{lemma}\label{gk}
(i) Suppose $m$ and $n$ are coprime odd integers, then
$$G_k(mn)=G_k(m)G_k(n).$$

(ii) Suppose $p^\alpha$ is the largest power of $p$ dividing $k$.
(If $k=0$ then set $\alpha=\infty$.)
Then for $\beta\geq 1$
\begin{equation*}
G_k(p^\beta)=\left\{\begin{aligned}
&0, &\beta\leq \alpha \textrm{ is odd},\\
&\phi(p^\beta), &\beta\leq \alpha \textrm{ is even},\\
&(\frac{kp^{-\alpha}}{p})p^\alpha\sqrt{p},  &\beta=\alpha+1 \textrm{ is odd},\\
&-p^\alpha, &\beta=\alpha+1 \textrm{ is even},\\
&0, &\beta\geq \alpha+2.
\end{aligned}
\right.
\end{equation*}
\end{lemma}

\begin{proof}
  See \cite[Lemma 2.3]{SD00}.
\end{proof}

For a Schwarz class function $F$ we define
\begin{equation} \label{eqn:tildeF}
  \widetilde F(\xi) 
  =\int_{-\infty}^{\infty}(\cos (2\pi \xi x)+\sin (2\pi \xi x))F(x) \dd x.
\end{equation}

\begin{lemma}[Poisson summation formula]\label{lm:psf}
Let $F$ be a nonnegative, smooth function supported in $(1,2)$.
For any odd integer $n$,
\begin{equation*}
S_M((\frac{d}{n});F)=\frac{1}{2n}(\frac{2}{n})\sum_{\substack{\alpha\leq Y\\ (\alpha,2n)=1}}\frac{\mu(\alpha)}{\alpha^2}\sum_{k=-\infty}^{\infty}(-1)^kG_k(n)\widetilde F(\frac{kX}{2\alpha^2 n}).
\end{equation*}
\end{lemma}

\begin{proof}
  See \cite[Lemma 2.6]{SD00}.
\end{proof}

\section{Proof of Theorem \ref{thm2}} \label{sec:first_moment}

By Lemma \ref{lemma:AFE}, we have
\begin{equation} \label{eqn:sum2S}
  \frac{1}{X} \sum_{2\nmid d}\nolimits^\flat L(\frac{1}{2},\phi\otimes \chi_{8d})\chi_{ 8d}(l)\Phi(\frac{d}{X})
  =
  S(\chi_{8d}(l)B(d);\Phi) + \prod_{i=1}^3\frac{
  \Gamma(\frac{\frac{1}{2}+\gamma_i}{2})} {\Gamma(\frac{\frac{1}{2}-\gamma_i}{2})}
  \overline{S(\chi_{8d}(l)B(d);\Phi)},
\end{equation}
where
\begin{equation*}
  B(d)=\sum_{n=1}^\infty \chi_{8d}(n)\frac{A(n,1)}{\sqrt{n}}V(n(\frac{\pi}{8d})^{\frac{3}{2}}).
\end{equation*}
We first consider the main contribution in $S(\chi_{8d}(l)B(d);\Phi)$, that is,
\begin{equation*}
  S_M(\chi_{8d}(l)B(d);\Phi)=
  \sum_{n=1}^\infty \frac{A(n,1)}{\sqrt{n}}S_M(\chi_{8d}(ln);\Phi_{n}),
\end{equation*}
where $\Phi_{y}(t)=\Phi(t)V(y(\frac{\pi}{8Xt})^\frac{3}{2})$.

By using Lemma \ref{lm:psf}, we obtain
\begin{equation*}
  S_M(\chi_{8d}(ln);\Phi_{n})=\frac{1}{2ln}(\frac{16}{ln})\sum_{\substack{\alpha \leq Y \\ (\alpha,2ln)=1}}\frac{\mu(\alpha)}{\alpha^2}\sum_{k=-\infty}^{\infty}(-1)^kG_k(ln)\widetilde \Phi_{n}(\frac{kX}{2\alpha^2 ln}).
\end{equation*}
Hence we deduce that
\begin{equation} \label{eqn:S(M)=P+R}
  S_M(\chi_{8d}(l)B(d);\Phi)=P(l)+R(l),
\end{equation}
where $P(l)$ are terms from $k=0$ and $R(l)$ are terms include all the nonzero terms $k$.
Thus
\begin{equation}\label{eqn:P}
P(l)=\frac{1}{2l}\sum_{n=1}^\infty \frac{A(n,1)}{n^\frac{3}{2}}(\frac{16}{ln})\sum_{\substack{\alpha \leq Y \\ (\alpha,2ln)=1}}\frac{\mu(\alpha)}{\alpha^2}G_0(ln)\widetilde \Phi_{n}(0),
\end{equation}
and
\begin{equation}\label{eqn:R}
R(l)=\frac{1}{2l}\sum_{n=1}^\infty \frac{A(n,1)}{n^\frac{3}{2}}(\frac{16}{ln})\sum_{\substack{\alpha \leq Y \\ (\alpha,2ln)=1}}\frac{\mu(\alpha)}{\alpha^2}\sum_{\substack{k=-\infty \\ k \neq 0}}^{\infty}(-1)^kG_k(ln)\widetilde \Phi_{n}(\frac{kX}{2\alpha^2 ln}).
\end{equation}

\subsection{The principal $P(l)$ contribution}\label{mainterm}
Note that $\widetilde \Phi_{n}(0)=\check{\Phi}_{n}(0)$ and that $G_0(ln)=\phi(ln)$ if $ln=\Box$ and $G_0(ln)=0$ otherwise.
Recall that $l=l_1l_2^2$ where $l_1$ and $l_2$ are odd, and $l_1$ is square-free.
The condition $ln=\Box$ is thus equivalent to $n=l_1m^2$ for some integer $m$.
Hence by \eqref{eqn:P} we have
\begin{equation}\label{eqn:P=}
\begin{split}
P(l)&
=\frac{1}{\zeta(2)\sqrt{l_1}}\sum_{\substack{m=1 \\ m \textrm{ odd}}}^\infty\frac{A(l_1m^2,1)}{m} (\prod_{p\mid 2lm}\frac{p}{p+1})\check{\Phi}_{l_1m^2}(0)\\
&\hskip 30pt +O(\frac{1}{Y\sqrt{l_1}}\sum_{\substack{m=1 \\ m \textrm{ odd}}}^\infty\frac{|A(l_1m^2,1)|}{m}|\check{\Phi}_{l_1m^2}(0)|).
\end{split}
\end{equation}
By Lemma \ref{lm2.2} and Theorem \ref{thm3}, together with some arguments as in  \cite[\S 2]{BV12}, we have
\begin{equation}\label{eqn:Perror}
\begin{split}
\sum_{\substack{m=1 \\ m \textrm{ odd}}}^\infty\frac{|A(l_1m^2,1)|}{m}|\check{\Phi}_{l_1m^2}(0)|&\ll \sum_{m^2\ll X^{\frac{3}{2}+\varepsilon}/l_1}\frac{|A(l_1m^2,1)|}{m}\\
&\ll \sum_{d\mid l_1^{\infty}} \sum_{\substack{m^2\ll X^{\frac{3}{2}+\varepsilon}/l_1d^2\\(m,l_1)=1}}\frac{|A(l_1d^2m^2,1)|}{dm}\\
&\ll \sum_{d\mid l_1^{\infty}} \frac{|A(l_1d^2,1)|}{d}\sum_{m^2\ll X^{\frac{3}{2}+\varepsilon}/l_1d^2}\frac{|A(m^2,1)|}{m}\\
&\ll l_1^{\theta_3+\varepsilon}\sum_{m\ll X^{\frac{3}{4}+\varepsilon}l_1^{-\frac{1}{2}+\varepsilon}}\frac{|A(m^2,1)|}{m}
\ll l_1^{\theta_3+\varepsilon}X^{\varepsilon}.
\end{split}
\end{equation}
By \eqref{eqn:P=} and \eqref{eqn:Perror}, we get
\begin{equation*}
  P(l) =\frac{1}{\zeta(2)\sqrt{l_1}}\sum_{\substack{m=1 \\ m \textrm{ odd}}}^\infty\frac{A(l_1m^2,1)}{m} (\prod_{p\mid 2lm}\frac{p}{p+1})\check{\Phi}_{l_1m^2}(0)
  + O(l_1^{\theta_3-\frac{1}{2}+\varepsilon}Y^{-1}).
\end{equation*}
For any $u>0$ we have
\begin{equation*}
\begin{split} \check{\Phi}_{l_1m^2}(0)&=\int_0^\infty \Phi(t) V(l_1m^2(\frac{\pi}{8Xt})^\frac{3}{2})\dd t\\
&=\frac{1}{2\pi i}\int_{(u)}\prod_{i=1}^3\frac{
\Gamma(\frac{s+\frac{1}{2}-\gamma_i}{2})}
{\Gamma(\frac{\frac{1}{2}-\gamma_i}{2})}(\frac{1}{l_1m^2}(\frac{8X}{\pi})^{\frac{3}{2}})^s (\int_0^\infty\Phi(t)t^{\frac{3s}{2}} dt)e^{s^2}\frac{\dd s}{s}\\
&=\frac{1}{2\pi i}\int_{(u)}\prod_{i=1}^3\frac{
\Gamma(\frac{s+\frac{1}{2}-\gamma_i}{2})} {\Gamma(\frac{\frac{1}{2}-\gamma_i}{2})} (\frac{1}{l_1m^2}(\frac{8X}{\pi})^{\frac{3}{2}})^s\check{\Phi}(\frac{3s}{2})e^{s^2} \frac{\dd s}{s}.
\end{split}
\end{equation*}
Thus
\begin{multline}\label{eqn:P==}
P(l)=
\frac{2}{3\zeta(2)\sqrt{l_1}}\frac{1}{2\pi i}\int_{(u)}\prod_{i=1}^3\frac{
\Gamma(\frac{s+\frac{1}{2}-\gamma_i}{2})}
{\Gamma(\frac{\frac{1}{2}-\gamma_i}{2})} (\frac{1}{l_1}(\frac{8X}{\pi})^{\frac{3}{2}})^s \check{\Phi}(\frac{3s}{2})\\
 \times\sum_{\substack{m=1 \\ m \textrm{ odd}}}^\infty\frac{A(l_1m^2,1)}{m^{1+2s}} (\prod_{p\mid lm} \frac{p}{p+1})e^{s^2}\frac{\dd s}{s}+O(l_1^{\theta_3-\frac{1}{2}+\varepsilon}Y^{-1}).
\end{multline}

Let $\alpha(p), \beta(p), \gamma(p)$ be the local parameters of $\phi$ at $p$,
and so $\alpha(p)\beta(p), \alpha(p)\gamma(p), \beta(p)\gamma(p)$ are the local parameters of $\widetilde\phi$ at $p$. Then we have
\begin{equation*}
\sum_{h=0}^\infty \frac{A(p^{h},1)}{p^{sh}}=(1-\alpha(p)p^{-s})^{-1}(1-\beta(p)p^{-s})^{-1}(1-\gamma(p)p^{-s})^{-1},
\end{equation*}
The local Euler factor of the symmetric square lift of $\phi$ is defined as
\begin{equation*}
\begin{split}
L_p(s,\sym^2 \phi)&=(1-\alpha(p)^2p^{-s})^{-1}(1-\beta(p)^2p^{-s})^{-1}(1-\gamma(p)^2p^{-s})^{-1}\\
&\hskip 30pt\times(1-\alpha(p)\beta(p)p^{-s})^{-1}(1-\alpha(p)\gamma(p)p^{-s})^{-1}(1-\beta(p)\gamma(p)p^{-s})^{-1}.
\end{split}
\end{equation*}
Let $S$ be a finite set of places of $\mathbb{Q}$. Define  $L^S(s,\sym^2 \phi)= \prod_{p\not\in S} L_p(s,\sym^2 \phi)$.  We have the following lemma.

\begin{lemma}\label{lemma:DS4P}
Suppose that $l=l_1l_2^2$ is as above. Then for $\Re (s)$ sufficiently large
\begin{equation}\label{eqn:DS-m}
\sum_{\substack{m=1 \\ m \textrm{ odd}}}^\infty\frac{A(l_1m^2,1)}{m^{s}}\prod_{p|lm} (\frac{p}{p+1})=\prod_{p\mid l}\frac{p}{p+1}G_\phi(s;l)L^{\{2\}}(s,\sym^2 \phi),
\end{equation}
where $G_\phi(s;l)=\prod_{\textrm{odd prime }p}G_{\phi,p}(s;l)$,
and
\begin{equation}\label{eqdefineta}
  G_{\phi,p}(s;l)=\left\{ \begin{aligned}
  & (A(p,1)+\frac{1}{p^s})
    (1-\frac{A(1,p)}{p^s}
     +\frac{A(p,1)}{p^{2s}}-\frac{1}{p^{3s}}),
     & \textrm{if } p\mid l_1,
     \\
  & (1+\frac{A(1,p)}{p^s})(1-\frac{A(1,p)}{p^s}
     +\frac{A(p,1)}{p^{2s}}-\frac{1}{p^{3s}}),
    &\textrm{if } p\mid l_2,p\nmid l_1,
    \\
  & (1-\frac{p^{2s}A(p,1)^2-
    p^{2s}A(1,p)-p^sA(1,p)^2+2p^sA(p,1)+1}
     {(p+1)p^{2s}(p^s+A(1,p))})
  \\
  & \hskip 30pt \times
(1+\frac{A(1,p)}{p^s})(1-\frac{A(1,p)}{p^s}
+\frac{A(p,1)}{p^{2s}}-\frac{1}{p^{3s}}),
&\textrm{if } p\nmid l.
\end{aligned}
\right.
\end{equation}
The right hand side of \eqref{eqn:DS-m} has analytic continuation to $\Re (s)>\frac{1}{2}$.
We have $|G_\phi(s;l)|\ll_{\sigma} l_1^{\theta_3+\varepsilon}$
if $\Re(s)=\sigma>1/2+\varepsilon$.
Moreover, there exists $c_\phi =3^a\times 5^b \times 7^2 \times 11^2$ with some $a,b\in\{1,2\}$
such that:
For any $l=c_\phi l'$ with $(2c_\phi,l')=1$,
$G_\phi(1;l)= 0$ if and only if there exists a prime $13\leq p\mid l_1$ with
$A(p,1)=-\frac{1}{p}$.
\end{lemma}

\begin{proof}
Expanding the Euler factors on the left, we get
\begin{equation}\label{eqn:DS2EP}
\begin{split} \sum_{\substack{m=1 \\ m \textrm{ odd}}}^\infty\frac{A(l_1m^2,1)}{m^{s}}
(\prod_{p\mid lm} \frac{p}{p+1})
&=\prod_{\substack{p \textrm{ prime}\\ p\mid l_1}}(\sum_{h=0}^\infty \frac{A(p^{2h+1},1)}{(p+1)p^{sh-1}})
\prod_{\substack{p \textrm{ prime}\\ p\mid l \\ p\nmid l_1}}(\sum_{h=0}^\infty \frac{A(p^{2h},1)}{(p+1)p^{sh-1}})\\
&\hskip 30pt \times \prod_{\substack{p \textrm{ prime}\\ p\nmid l}}(1+\sum_{h=1}^\infty \frac{A(p^{2h},1)}{(p+1)p^{sh-1}}).
\end{split}
\end{equation}
Recall from \cite{GB06} the relationship between the coefficients of $\phi$:
\begin{equation*}
A(m_1,1)A(1,m_2)=\sum_{d\mid (m_1,m_2)}A(\frac{m_1}{d},\frac{m_2}{d}),
\end{equation*}
\begin{equation*}
A(1,n)A(m_1,m_2)=\sum_{\substack{d_0d_1d_2=n \\ d_1\mid m_1 \\ d_2\mid m_2}}A(\frac{m_1d_2}{d_1},\frac{m_2d_0}{d_2}),
\end{equation*}
so we have
\begin{equation}\label{heckerealtion}
\begin{split}
A(p^{k+1},1)=A(p,1)A(p^{k},1)-A(p^{k-1},p),\\
A(p^{k-1},p)=A(p^{k-1},1)A(1,p)-A(p^{k-2},1),
\end{split}
\end{equation}
thus
\begin{equation}\label{eqn:sum-h1}
\sum_{h=0}^\infty \frac{A(p^{2h+1},1)}{(p+1)p^{sh-1}}=\frac{p}{p+1}\frac{1+p^sA(p,1)}{p^s+A(1,p)}\sum_{h=0}^\infty \frac{A(p^{2h},1)}{p^{sh}},
\end{equation}
and
\begin{equation}\label{eqn:sum-h2}
\begin{split}
1+\sum_{h=1}^\infty \frac{A(p^{2h},1)}{(p+1)p^{sh-1}}
&=\sum_{h=0}^\infty \frac{A(p^{2h},1)}{p^{sh}}-\frac{1}{p+1}\sum_{h=1}^\infty \frac{A(p^{2h},1)}{p^{sh}}\\
&=\sum_{h=0}^\infty \frac{A(p^{2h},1)}{p^{sh}}-\frac{1}{(p+1)p^s}\sum_{h=0}^\infty \frac{A(p^{2h+2},1)}{p^{sh}}\\
&=(1+\frac{A(1,p)}{(p+1)p^s})\sum_{h=0}^\infty \frac{A(p^{2h},1)}{p^{sh}}-\frac{1+p^sA(p,1)}{(p+1)p^{2s}}\sum_{h=0}^\infty \frac{A(p^{2h+1},1)}{p^{sh}}\\
&=(1-\frac{p^{2s}A(p,1)^2-p^{2s}A(1,p)-p^sA(1,p)^2+2p^sA(p,1)+1}{(p+1)p^{2s}(p^s+A(1,p))})\\
&\hskip 30pt\times\sum_{h=0}^\infty \frac{A(p^{2h},1)}{p^{sh}}.
\end{split}
\end{equation}
Recall the Euler factors of $L(s,\phi)$ we know
\begin{equation*}
A(p^h,1)=\sum_{a+b+c=h}\alpha(p)^a\beta(p)^b\gamma(p)^c.
\end{equation*}
Note that if the sum of three integers is an even integer, then there will be zero or two odd integers.  So we have
\begin{equation*}
\begin{split}
A(p^{2h},1)&=\sum_{a+b+c=2h}\alpha(p)^a\beta(p)^b\gamma(p)^c\\
&=\sum_{a+b+c=h}\alpha(p)^{2a}\beta(p)^{2b}\gamma(p)^{2c} \\
& \qquad +(\alpha(p)\beta(p)+\alpha(p)\gamma(p)+\beta(p)\gamma(p)) \sum_{a+b+c=h-1}\alpha(p)^{2a}\beta(p)^{2b}\gamma(p)^{2c},
\end{split}
\end{equation*}
thus
\begin{equation*}
\begin{split}
\sum_{h=0}^\infty \frac{A(p^{2h},1)}{p^{sh}}
&=(1+\frac{\alpha(p)\beta(p)+\alpha(p)\gamma(p)+\beta(p)\gamma(p)}{p^s})\sum_{h=0}^\infty \frac{\sum_{a+b+c=h}\alpha(p)^{2a}\beta(p)^{2b}\gamma(p)^{2c}}{p^{sh}}\\
&=\frac{p^s+\alpha(p)\beta(p)+\alpha(p)\gamma(p)+\beta(p)\gamma(p)} {p^s} \\
& \qquad \cdot (1-\alpha(p)^2p^{-s})^{-1}(1-\beta(p)^2p^{-s})^{-1}(1-\gamma(p)^2p^{-s})^{-1}.
\end{split}
\end{equation*}
From
\begin{equation*}
\begin{split}
(1-\alpha(p)^2p^{-s})^{-1}(1-\beta(p)^2p^{-s})^{-1}(1-\gamma(p)^2p^{-s})^{-1}&=L_p(s,\phi\otimes \phi)/L_p(s,\widetilde \phi)^2\\
&=L_p(s,\sym^2 \phi)/L_p(s,\widetilde \phi),
\end{split}
\end{equation*}
and
\begin{equation*}
\alpha(p)\beta(p)+\alpha(p)\gamma(p)
+\beta(p)\gamma(p)=A(1,p),
\end{equation*}
we obtain
\begin{equation} \label{eqn:sum-2h}
\begin{split}
\sum_{h=0}^\infty \frac{A(p^{2h},1)}{p^{sh}}&=(1+\frac{A(1,p)}{p^s})L_p(s,\widetilde \phi)^{-1}L_p(s,\sym^2 \phi)\\
&=(1+\frac{A(1,p)}{p^s})(1-\frac{A(1,p)}{p^s}
+\frac{A(p,1)}{p^{2s}}-\frac{1}{p^{3s}})
L_p(s,\sym^2 \phi).
\end{split}
\end{equation}
By \eqref{eqn:DS2EP}--\eqref{eqn:sum-2h}, we prove \eqref{eqn:DS-m}.

Recall that we have $\theta_3\leq\frac{5}{14}$, thus when $\Re (s)=\sigma>\frac{1}{2}$ we have
\begin{equation*}
\begin{split}
\log \prod_
{\substack{2< p\leq Z\\p\nmid l_1}}
G_{\phi,p}(s;l)
&\ll\sum_{p\leq Z}|A(1,p)|^2p^{-2\sigma}+|A(p,1)|p^{-2\sigma}
+|A(p,1)|^2 p^{-1-\sigma} + |A(p,1)|p^{-1-\sigma}\\
&\ll\sum_{n\leq Z}|A(1,n)|^2n^{-2\sigma}+|A(1,n)|n^{-2\sigma}
+|A(1,n)|^2 n^{-1-\sigma} + |A(n,1)|n^{-1-\sigma},
\end{split}
\end{equation*}
so $\prod_{{\substack{2< p\leq Z\\p\nmid l_1}}}|G_{\phi,p}(s;l)|\ll_{\sigma} 1$ for $\Re (s)>\frac{1}{2}$.
For fixed $l$ there only finite $p\mid l_1$,
thus $G_{\phi}(s;l)$ converges for $\Re (s)>\frac{1}{2}$.
Moreover, we have $|G_\phi(s;l)|\ll_{\phi,\sigma} l_1^{\theta_3+\varepsilon}$
if $\Re(s)=\sigma>1/2+\varepsilon$ by using the known Ramanujan bounds to the local factors at primes dividing $l_1$.

Finally, we will prove the last claim.
It is known that $L_p(1,\widetilde{\phi})^{-1} \neq 0$, since we have $L_p(1,\widetilde{\phi})^{-1} \sum_{h\geq 0} \frac{A(1,p^h)}{p^h} = 1$.
From $|A(p,1)|\leq 3p^{\frac{5}{14}}$ we know that for $p\geq 13$,
\[ |A(1,p)|/p < 1 \]
and
\begin{equation*}
|\frac{p^{2}A(p,1)^2-p^{2}A(1,p)-pA(1,p)^2
+2pA(p,1)+1}{(p+1)p^{2}(p+A(1,p))}|<1.
\end{equation*}
So $G_{\phi,p}(1;l)\neq 0$ for $13\leq p\nmid l$. Note that  we have
\begin{equation*}
\begin{split}
\log \prod_{\substack{13\leq p\leq Z\\p\nmid l}}|G_{\phi,p}(1;l)|&\gg -\sum_{13\leq p\leq Z}\Big(|A(1,p)|^2 p^{-2}+|A(p,1)|
p^{-2}\Big)\\
&\gg -\sum_{n\leq Z}|A(1,n)|^2n^{-2}-\sum_{n\leq Z}|A(1,n)|n^{-2},
\end{split}
\end{equation*}
which proves $\prod_{13\leq p \nmid l}G_{\phi,p}(1;l)\neq 0$.

For $p\mid l_1$, $G_{\phi,p}(1,l)=0$ if and only if $A(p,1)=-\frac{1}{p}$;
and for $p\mid l_2$, $p\nmid l_1$, we have $G_{\phi,p}(1,l)=0$ if and only if $A(p,1)=-p$ which may happen only if $p=3$ or $5$.
Thus we know there exists $c_\phi =3^a\times 5^b \times 7^2 \times 11^2$ with some $a,b\in\{1,2\}$
such that for any $l=c_\phi l'$ with $(2c_\phi ,l')=1$,
$\prod_{3\leq p\leq 11}G_{\phi,p}(1;l)\neq 0$.
For such $l$, we have $G_{\phi}(1;l)= 0$ if and only if there is one prime $13\leq p\mid l_1$ such that
$A(p,1)=-\frac{1}{p}$.
This completes the proof of the lemma.
\end{proof}


We denote $G_{\phi}(l)=G_{\phi}(1;l)$.
By \eqref{eqn:P==} and Lemma \ref{lemma:DS4P}, we have
\begin{equation}\label{eqn:P===}
P(l)=\frac{2}{3}\frac{1}{\zeta(2)\sqrt{l_1}}\prod_{p\mid l}\frac{p}{p+1}I(l) +O(l_1^{\theta_3-\frac{1}{2}+\varepsilon}Y^{-1}),
\end{equation}
where
\begin{equation*}
\begin{split}
I(l)=&\frac{1}{2\pi i}\int_{(u)}\prod_{i=1}^3\frac{
\Gamma(\frac{s+\frac{1}{2}-\gamma_i}{2})}
{\Gamma(\frac{\frac{1}{2}-\gamma_i}{2})}
(\frac{1}{l_1}
(\frac{8X}{\pi})^{\frac{3}{2}})^s\check{\Phi}
(\frac{3s}{2})e^{s^2} \\
&\hskip 30pt \times G_{\phi}(1+2s;l)
L^{\{2\}}(1+2s,\sym^2 \phi)\frac{\dd s}{s}.
\end{split}
\end{equation*}

We move the line of integration to the $\Re (s)=-u$ line with $u=\min\{\frac{1}{4},\frac{1}{2}-\theta_3\}-\varepsilon$.
From \cite{GB06} we know there is a pole of the integrand at $s=0$ and we shall evaluate the residue of this pole shortly.
We now bound the integral on the $-u$ line.
From \cite{IK04} we know that on this line we have
\begin{equation*}
|L(1+2s,\sym^2 \phi)|\ll \prod_{i=1}^3(|s|+|\gamma_i|+3)^6, \quad -\frac{1}{4}<\Re s<0,
\end{equation*}
\begin{equation*}
|L_2(1+2s,\sym^2 \phi)|\gg (1-2^{-1-2\Re(s)+2\theta_3})^6,
\end{equation*}
and on the $-u$ line we have $|G_{\phi}(1+2s,l)|\ll l_1^{\theta_3+\varepsilon}$.
Hence the integral on the line is
\begin{equation}\label{eqn:Perror2}
\begin{split}
&\ll \frac{l_1^{u+\theta_3+\varepsilon}}{X^{\frac{3u}{2}-\varepsilon}}\int_{(-u)}
\frac{\prod_{i=1}^3(|s|+|\gamma_i|+3)^6}
{(1-2^{-1-2\Re(s)+2\theta_3})^6}|
\check{\Phi}(\frac{3s}{2})||e^{s^2}|
|\prod_{i=1}^3\Gamma(\frac
{s+\frac{3}{2}-\gamma_i}{2})|\frac{|\dd s|}{|s|}\\
&\ll \frac{l_1^{u+\theta_3+\varepsilon}}{X^{\frac{3u}{2}-\varepsilon}}.
\end{split}
\end{equation}
When $\phi$ is not self-dual, $L(s,\sym^2 \phi)$ has no pole or zero point at $s=1$, we evaluate residues of pole at $s=0$ are
$\hat{\Phi}(0)G_{\phi}(l)L^{\{2\}}(1,\sym^2 \phi)$.
When $\phi$ is self-dual, we know $L(s,\sym^2 \phi)$ has a simple pole at $s=1$, so we have the Laurent series expansions
\begin{equation*}
\begin{gathered}
\prod_{i=1}^3\frac{
\Gamma(\frac{s+\frac{1}{2}-\gamma_i}{2})}
{\Gamma(\frac{\frac{1}{2}-\gamma_i}{2})}
=1+as+\dots;\\
(\frac{1}{l_1}(\frac{8X}{\pi})
^{\frac{3}{2}})^s=1+\frac{3}{2}
\log(\frac{8X}{l_1^\frac{2}{3}\pi})s+\dots;\\
G_\phi(1+2s;l)=G_\phi(l)
+2G^{'}_\phi(1;l)s+\dots;\\
 L^{\{2\}}(1+2s,\sym^2 \phi) =\lim_{s_1\to 0} s_1 L^{\{2\}}(1+2s_1,\sym^2 \phi)  \frac{1}{s}+c_1+c_2s+\dots;
\end{gathered}
\end{equation*}
and $\check{\Phi}(\frac{3s}{2})e^{s^2}=\check{\Phi}(0)+\frac{3}{2}\check{\Phi}^{'}(0)s+\dots$.
It follows that the residues may be written as
\begin{equation}\label{cphi}
\frac{3}{2}\lim_{s_1\to 0}s_1L^{\{2\}}(1+2s_1,\sym^2 \phi)\check{\Phi}(0)
\Big(G_\phi(l)
\log \frac{X}{l_1 ^{\frac{2}{3}}}+C_\phi(l)\Big)
\end{equation}
where
\begin{equation}\label{Cphi}
  C_\phi(l)=G_\phi(l)(\frac{2}{3}a-\log \pi+\frac{\check{\Phi}^{'}(0)}
   {\check{\Phi}(0)})
  +\frac{4}{3}G^{'}_\phi(1;l).
\end{equation}

By  \eqref{eqn:P===},  \eqref{eqn:Perror2}, and \eqref{cphi}, we conclude that if $\phi$ is not self-dual, then
\begin{multline} \label{eqn:P=1}
P(l)=
\frac{2\check{\Phi}(0)}{3\zeta(2)\sqrt{l_1}}\prod_{p\mid l}\frac{p}{p+1}G_\phi(l)L^{\{2\}}(1,\sym^2 \phi)\\
+O(\min\{l_1^{-\frac{1}{4}+\theta_3+\varepsilon} X^{-\frac{3}{8}+\varepsilon},l_1^{\frac{1}{2}+\varepsilon} X^{\frac{3}{2}\theta_3-\frac{3}{4}+\varepsilon}\}) +O(l_1^{\theta_3-\frac{1}{2}+\varepsilon}Y^{-1});
\end{multline}
and if $\phi$ is self-dual, then
\begin{multline} \label{eqn:P=2}
P(l)=\lim_{s_1\to 0}s_1L^{\{2\}}(1+2s_1,\sym^2 \phi)\frac{\check{\Phi}(0)}{\zeta(2)\sqrt{l_1}}
\prod_{p\mid l}\frac{p}{p+1}\Big(G_\phi(l)
\log \frac{X}{l_1 ^{\frac{2}{3}}}+C_\phi(l)\Big)\\
+O(\min\{l_1^{-\frac{1}{4}+\theta_3+\varepsilon}
X^{-\frac{3}{8}+\varepsilon},
l_1^{\frac{1}{2}+\varepsilon}
X^{\frac{3}{2}\theta_3-\frac{3}{4}+\varepsilon}\})
+O(l_1^{\theta_3-\frac{1}{2}+\varepsilon}Y^{-1}).
\end{multline}

\subsection{The contribution of the remainder terms $R(l)$}\label{remainder1}

By using inverse Mellin transform, we have
\begin{equation*}
\sum_{n=1}^\infty a_ng(n)=\frac{1}{2\pi i}\int_{(c)}\sum_{n=1}^\infty \frac{a_n}{n^w}(\int_0^\infty g(t)t^{w-1}\dd t) \dd w.
\end{equation*}
By \eqref{eqn:R}, we may recast the expression for $R(l)$ as
\begin{equation}\label{eqn:R=}
R(l)=\frac{1}{2l}\sum_{\substack{\alpha \leq Y \\ (\alpha,2l)=1}}\frac{\mu(\alpha)}{\alpha^2}\sum_{\substack{k=-\infty \\ k \neq 0}}^{\infty}\frac{(-1)^k}{2\pi i}\int_{(c)}\sum_{\substack{n=1\\ (n,2\alpha)=1}} ^\infty \frac{A(n,1)}{n^{\frac{3}{2}+w}}G_{4k}(ln)\phi(\frac{kX}{2\alpha^2 l},w)\dd w
\end{equation}
for any $c>0$, where
\begin{equation} \label{eqn:phi(xi&w)}
\phi(\xi,w)=\int_0^\infty \widetilde \Phi_{t}(\frac{\xi}{t})t^{w-1}\dd t
\end{equation}
 with
 \begin{equation}\label{eqn:Phi_t}
 \Phi_{y}(t)=\Phi(t)V(y(\frac{\pi}{8Xt})^\frac{3}{2}).
\end{equation}

To estimate $R(l)$, we will need the following results for the sum over $n$ (Lemma \ref{lemma:DS-n}) and for $\phi(\xi,w)$ (Lemma \ref{lemma:phi}).

\begin{lemma}\label{lemma:DS-n}
Write $4k=k_1k_2^2$ where $k_1$ is a fundamental discriminant (possibly $k_1=1$ is the trivial character), and $k_2$ is positive.
In the region $\Re (s)>1$ we have
\begin{equation*}
\begin{split}
\sum_{\substack{n=1\\ (n,2\alpha)=1}}^\infty \frac{A(n,1)}{n^{s}}\frac{G_{4k}(ln)}{\sqrt{n}}
&=L(s,\phi\otimes \chi_{k_1})H(s;\phi,k,l,\alpha),
\end{split}
\end{equation*}
where $H(s;\phi,k,l,\alpha)$ has an analytic continuation to $\Re(s)>1/2$.
For $\Re (s)>\frac{1}{2}+\varepsilon$ we have
\begin{equation}\label{eqn:H<<}
  |H(s;\phi,k,l,\alpha)|\ll |k|^\varepsilon \alpha^\varepsilon l^{\frac{1}{2}+\varepsilon}
(l,k_2^2)^{\frac{1}{2}}
\sum_{\substack{h_0\mid \frac{k\cdot \textrm{rad } k}{(l,k)}\\p\mid h_0\Rightarrow p\mid (l,k_2)}}
\sum_{\substack{h\mid k\\(h,l)=1}}|A(h_0h,1)|,
\end{equation}
where $\textrm{rad } n=\prod_{p\mid n}p$ is the radical of $n$.
\end{lemma}

\begin{proof}
  By the multiplicativity of $G_{4k}(n)$, we have the following Euler product expansion
\begin{equation*}
\begin{split}
\sum_{\substack{n=1\\ (n,2\alpha)=1}}^\infty \frac{A(n,1)}{n^{s}}\frac{G_{4k}(ln)}{\sqrt{n}}
&=L(s,\phi\otimes \chi_{k_1})\prod_p H_{p}(s;\phi,k,l,\alpha)\\
&=L(s,\phi\otimes \chi_{k_1})H(s;\phi,k,l,\alpha),
\end{split}
\end{equation*}
where $H_{p}$ is defined as follows:
\begin{equation*}
H_{p}(s;\phi,k,l,\alpha)=\left\{\begin{aligned}
&(1-A(p,1)(\frac{k_1}{p})p^{-s}+A(1,p)(\frac{k_1}{p})^2 p^{-2s}
-(\frac{k_1}{p})p^{-3s}), &p\mid 2\alpha,\\
&(1-A(p,1)(\frac{k_1}{p})p^{-s}+A(1,p)(\frac{k_1}{p})^2 p^{-2s}
-(\frac{k_1}{p})p^{-3s})\\
&\hskip 30pt \times \sum_{r=0}^{\infty}\frac{A(p^r,1)}{p^{rs}}\frac{G_k(p^{r+ord_p (l)})}{p^{\frac{r}{2}}}, &p\nmid 2\alpha.
\end{aligned}
\right.
\end{equation*}
We see that for a generic $p\nmid 2\alpha kl$,
from Lemma \ref{gk} we have $G_k(p^{r+ord_p(l)})=0$ for $r\geq 2$,
so
\begin{equation*}
\begin{split}
H_{p}(s;\phi,k,l,\alpha)
&=(1-A(p,1)(\frac{k_1}{p})p^{-s}+A(1,p)(\frac{k_1}{p})^2 p^{-2s}-(\frac{k_1}{p})p^{-3s}) (1+(\frac{k_1}{p}) \frac{A(p,1)}{p^s}) \\
&=1+(A(1,p)-A(p,1)^2)p^{-2s}+(A(p,1)A(1,p)-1)(\frac{k_1}{p})p^{-3s}-A(p,1)p^{-4s}\\
&=1-A(p^2,1)p^{-2s}+A(p,p)(\frac{k_1}{p})p^{-3s}-A(p,1)p^{-4s}.
\end{split}
\end{equation*}
Note that for $\Re (s)=\sigma>\frac{1}{2}+\varepsilon$
\begin{equation*}
\begin{split}
\log \prod_{\substack{p\leq Z\\p\nmid 2\alpha kl}}H_{p}(s;\phi,k,l,\alpha)&\ll_\varepsilon \sum_{p\leq Z}|A(p^2,1)|p^{-2\sigma}+|A(p,p)|p^{-3\sigma}+|A(p,1)|p^{-4\sigma}\\
&\ll_\varepsilon \sum_{p\leq Z}(|A(p,1)^2|+|A(p,1)|)p^{-2\sigma}\\
&\ll_\varepsilon \sum_{n\leq Z}(|A(n,1)^2|+|A(n,1)|)n^{-2\sigma}\\
&\ll_\varepsilon 1.
\end{split}
\end{equation*}
This shows that $H(s;\phi,k,l,\alpha)$ is holomorphic in $\Re (s)=\sigma>\frac{1}{2}+\varepsilon$.

It remains to prove the bound \eqref{eqn:H<<}.
From our evaluation of $H_{p}(s;\phi,k,l,\alpha)$ for $p\nmid 2kl\alpha$ we see that for $\Re (s)>\frac{1}{2}+\varepsilon$,
\begin{equation*}
|H(s;\phi,k,l,\alpha)|\ll (|k|l\alpha)^\varepsilon \prod_{\substack{p\mid kl\\ p\nmid 2\alpha}}|H_{p}(s;\phi,k,l,\alpha)|.
\end{equation*}
Suppose now that $p^a\parallel k$ and $p^b\parallel l$ with $a+b\geq1$ and $p\nmid 2\alpha$.
By Lemma \ref{gk}, we may suppose that $b\leq a+1$, since otherwise we have  $H_{p}(s;\phi,k,l,\alpha)=0$.
Note that $|G_k(p^{j})|\leq p^{j}$ for $0\leq j\leq a$,
and consider the cases of $a$ even and $b=a+1$; and $a$ odd and $b=a$, or $b=a+1$, then we have
\begin{equation*}
\begin{split}
|H_{p}(s;\phi,k,l,\alpha)| & \leq (1+3p^{-\frac{1}{7}}
+3p^{-\frac{9}{14}}+p^{-\frac{3}{2}})
p^{\min\{b,[\frac{a}{2}]+\frac{b}{2}\}}
\sum_{r=0}^{a+1-b}\frac{|A(p^r,1)|}
{p^{r(\sigma-\frac{1}{2})}} .
\end{split}
\end{equation*}
By the Hecke relation \eqref{heckerealtion}, $|G_k(p^{a+1})|\leq p^{a+\frac{1}{2}}$, and the bounds toward to Ramanujan conjecture when $b=0$, we have
\begin{equation*}
\begin{split}
|H_{p}(s;\phi,k,l,\alpha)|
& \leq (2+ 3p^{-\frac{1}{7}}
+3p^{-\frac{9}{14}}+p^{-\frac{3}{2}} )
p^{\min\{b,[\frac{a}{2}]+\frac{b}{2}\}}
\sum_{r=0}^{a}|A(p^r,1)|.
\end{split}
\end{equation*}
For $a=1$ from $|G_k(p^{2})|=p$ we also have $|H_{p}(s;\phi,k,l,\alpha)|\leq 2 p^{\frac{b}{2}}$ to handle with $p\mid (k_1,l),p\nmid k_2$.
Hence we get
\begin{equation*}
|H(s;\phi,k,l,\alpha)|\ll |k|^\varepsilon \alpha^\varepsilon l^{\frac{1}{2}+\varepsilon}
(l,k_2^2)^{\frac{1}{2}}
\sum_{\substack{h_0\mid \frac{k\cdot \textrm{rad } k}{(l,k)}\\p\mid h_0\Rightarrow p\mid (l,k_2)}}
\sum_{\substack{h\mid k\\(h,l)=1}}|A(h_0h,1)|,
\end{equation*}
from which we finish the proof.
\end{proof}

\begin{lemma}\label{lemma:phi}
We have
\begin{equation*}
\begin{split}
\phi(\xi,w)&=\frac{1}{2\pi i}\int_{(u)}\bigg(\cos\big(\frac{\pi}{2} (s-w)\big)+\mathrm{sgn}(\xi)\sin\big(\frac{\pi}{2} (s-w)\big)\bigg)(\frac{8X}{\pi})^{\frac{3s}{2}}(2\pi|\xi|)^{w-s} \Gamma(s-w)\\&\hskip 30pt \times\prod_{i=1}^3\frac{
\Gamma(\frac{s+\frac{1}{2}-\gamma_i}{2})} {\Gamma(\frac{\frac{1}{2}-\gamma_i}{2})}
\check{\Phi}(\frac{s}{2}+w)\frac{e^{s^2}}{s}\dd s.
\end{split}
\end{equation*}
\end{lemma}

\begin{proof}
By  \eqref{eqn:tildeF} and \eqref{eqn:phi(xi&w)}, we have
\begin{equation*}
  \phi(\xi,w)=\int_0^\infty \left(
  \int_0^\infty \Phi(y)V(t(\frac{\pi}{8Xy})^\frac{3}{2}) (\cos(2\pi y\frac{\xi}{t})+\sin(2\pi y\frac{\xi}{t}))\dd y \right)
  t^{w-1}\dd t.
\end{equation*}
In the inner integral over $y$, we make the substitution $z=|\xi|y/t$, so that this integral becomes
\begin{equation*}
\frac{t}{|\xi|}\int_0^\infty\Phi_{t}(\frac{tz}{|\xi|})(\cos(2\pi z)+\mathrm{sgn}(\xi)\sin(2\pi z))\dd z.
\end{equation*}
We use this above, and interchange the integrals over $z$ and $t$.
Thus
\begin{equation*}
\phi(\xi,w)=\frac{1}{|\xi|}\int_0^\infty(\int_0^\infty\Phi_{t}(\frac{tz}{|\xi|})t^{w}dt)(\cos(2\pi z)+\mathrm{sgn}(\xi)\sin(2\pi z))\dd z.
\end{equation*}
From the definitions of $\Phi_{t}$ \eqref{eqn:Phi_t} and $V$  \eqref{eqn:V}, the inner integral is
\begin{equation*}
\begin{split}
\int_0^\infty V(t(\frac{\pi|\xi|}{8Xtz})
^\frac{3}{2})\Phi(\frac{tz}{|\xi|})t^{w}\dd t
&=\frac{1}{2\pi i}\int_0^\infty \int_{(u)}\prod_{i=1}^3\frac{
\Gamma(\frac{s+\frac{1}{2}-\gamma_i}{2})} {\Gamma(\frac{\frac{1}{2}-\gamma_i}{2})} (\frac{8Xz}{|\xi|\pi})^{\frac{3s}{2}} t^{\frac{s}{2}+w}\Phi(\frac{tz}{|\xi|})
e^{s^2}\frac{\dd s}{s}\dd t\\
&=\frac{1}{2\pi i}\int_{(u)}\prod_{i=1}^3\frac{
\Gamma(\frac{s+\frac{1}{2}-\gamma_i}{2})} {\Gamma(\frac{\frac{1}{2}-\gamma_i}{2})} (\frac{8Xz}{|\xi|\pi})^{\frac{3s}{2}}
\left( \int_0^\infty t^{\frac{s}{2}+w}\Phi(\frac{tz}{|\xi|})\dd t \right)
e^{s^2}\frac{\dd s}{s}\\
&=\frac{1}{2\pi i}\int_{(u)} \prod_{i=1}^3\frac{
\Gamma(\frac{s+\frac{1}{2}-\gamma_i}{2})} {\Gamma(\frac{\frac{1}{2}-\gamma_i}{2})} \check{\Phi}(\frac{s}{2}+w)(\frac{8Xz}
{|\xi|\pi})^{\frac{3s}{2}} (\frac{|\xi|}{z})^{\frac{s}{2}+w+1}
e^{s^2}\frac{\dd s}{s}.
\end{split}
\end{equation*}
Here in the last equality, we have made a change of variable from $tz/|\xi|$ to $t$ and used the definition of $\check{\Phi}$.
Thus
\begin{equation*}
\begin{split}
\phi(\xi,w)&=\frac{1}{2\pi i}\int_0^\infty\int_{(u)}(\cos(2\pi z)+\mathrm{sgn}(\xi)\sin(2\pi z))(\frac{8X}{\pi})^{\frac{3s}{2}}
|\xi|^{w-s}z^{s-w-1}\\
&\hskip 30pt\times \frac{e^{s^2}}{s}
\prod_{i=1}^3\frac{\Gamma(
\frac{s+\frac{1}{2}-\gamma_i}{2})} {\Gamma(\frac{\frac{1}{2}-\gamma_i}{2})}
\check{\Phi}(\frac{s}{2}+w)\dd s  \dd z.
\end{split}
\end{equation*}
Interchange the integrals over $s$ and $z$, employing the expressions for the Fourier sine and cosine transforms of $z^{s-w-1}$, we obtain the lemma.
\end{proof}

From \cite{GB06} we know $L$-functions for Hecke--Maass forms are all entire.
By Lemmas \ref{lemma:DS-n} and \ref{lemma:phi} and moving the lines of $w$ and $s$ such that $\Re (w-s)=-\frac{5}{4}-2\varepsilon$ and $\Re (w)=-\frac{1}{2}+2\varepsilon$ (so $\Re (s)=\frac{3}{4}+4\varepsilon$) in \eqref{eqn:R=}, we obtain
\begin{equation*}
\begin{split}
R(l)&=\frac{1}{4l\pi i}\sum_{\substack{\alpha \leq Y \\ (\alpha,2l)=1}}\frac{\mu(\alpha)}{\alpha^2}
\sum_{\substack{k=-\infty \\ k \neq 0}}^{\infty}\frac{(-1)^k}{2\pi i}
\int_{(-\frac{1}{2}+2\varepsilon)}
\int_{(\frac{3}{4}+4\varepsilon)}
L(1+w,\phi\otimes \chi_{k_1}) H(1+w;\phi,k,l,\alpha)
\\
&\hskip 30pt \times (\frac{k}{\alpha^2 l})^{w-s}X^{\frac{s}{2}+w}\pi^{w-
\frac{5s}{2}}8^{\frac{3s}{2}}\bigg(
\cos\big(\frac{\pi}{2} (s-w)\big)+\mathrm{sgn}(k)\sin\big(
\frac{\pi}{2} (s-w)\big)\bigg)\\
&\hskip 60pt \times
\prod_{i=1}^3\frac{\Gamma(\frac{s+
\frac{1}{2}-\gamma_i}{2})}
{\Gamma(\frac{\frac{1}{2}-\gamma_i}{2})} \check{\Phi }(\frac{s}{2}+w)
\Gamma(s-w)\frac{e^{s^2}}{s}\dd s \dd w \\
&\ll \frac{l^{\frac{3}{4}+3\varepsilon}}
{X^{\frac{1}{8}-4\varepsilon}}
\sum_{\alpha\leq Y} \alpha^{1/2+\varepsilon}
\int_{(\frac{3}{4}+4\varepsilon)}
\int_{(-\frac{1}{2}+2\varepsilon)}
\sum_{k_2=1}^{\infty}
\sum_{\substack{Z\geq 1 \\ \mathrm{dyadic}}}
\sum_{Z\leq k_1\leq 2Z}|L(1+w,\phi\otimes \chi_{k_1})|
(l,k_2^2)^{\frac{1}{2}}(Zk_2^{2})^{-\frac{5}{4}-\varepsilon}
  \\
& \hskip 30pt \times
\sum_{h\mid k_1k_2^{2}}|A(h,1)|
|\check{\Phi }(\frac{s}{2}+w)|(1+|s-w|)
^{\frac{3}{4}+2\varepsilon}\prod_{i=1}^3
|\Gamma(\frac{s+\frac{1}{2}-\gamma_i}{2})| |\frac{e^{s^{2}}}{s}| |\dd w \dd s|.
\end{split}
\end{equation*}
Here we have used the Stirling's formula to estimate $\Gamma(s-w)$.
We have
\begin{equation*}
\sum_{\substack{h_0\mid \frac{k\cdot \textrm{rad } k}{(l,k)}\\p\mid h_0\Rightarrow p\mid (l,k_2)}}
\sum_{\substack{h\mid k\\(h,l)=1}}|A(h_0h,1)|\leq
\sum_{h_1\mid k_1}
\sum_{\substack{h_2\mid \frac{k_2^{3}}{(l,k_2^2)}\\p\mid h_2\Rightarrow p\mid l}}
\sum_{\substack{h_3\mid k_2^{2}\\(h_3,l)=1}}|A(h_1h_2h_3,1)|.
\end{equation*}
By using the approximate functional equations (cf. Lemma \ref{lemma:AFE}) and the large sieve estimate for quadratic characters in \cite{HB95}, for $\Re(w)=-1/2+2\varepsilon$,
and in a similar way as \eqref{eqn:Perror},
we get
\begin{equation*}
\begin{split}
\sum_{Z\leq k_1\leq 2Z}&|L(1+w,\phi\otimes \chi_{k_1})| \sum_{h_1\mid k_1}|A(h_1h_2h_3,1)| \\
  &\ll (\sum_{Z\leq k_1\leq 2Z}|L(1+w,\phi\otimes \chi_{k_1})|^2)^\frac{1}{2}
  (2Z\sum_{k_1\leq 2Z}\sum_{h_1\mid k_1}\frac{|A(h_1h_2h_3,1)|^2}{k_1})^\frac{1}{2}
  \\
  &\ll Z^{\frac{5}{4}+\frac{\varepsilon}{2}}
  (3+|w|+\sum_{i=1}^3|\gamma_i|)
  ^{\frac{3}{4}+\varepsilon}
    (h_2h_3)^{\theta_3+\varepsilon}.
\end{split}
\end{equation*}
Hence we have
\begin{equation} \label{eqn:R<<}
\begin{split}
R(l)
&
\ll \frac{l^{\frac{3}{4}+3\varepsilon}Y^{\frac{3}{2}+2\varepsilon}}{X^{\frac{1}{8}-4\varepsilon}}
\int_{(\frac{3}{4}+4\varepsilon)}
\int_{(-\frac{1}{2}+2\varepsilon)}
\sum_{k_2=1}^{\infty}
\sum_{\substack{h_2\mid \frac{k_2^{3}}{(l,k_2^2)}\\p\mid h_2\Rightarrow p\mid l}}
\sum_{\substack{h_3\mid k_2^{2}\\(h_3,l)=1}}
\frac{(l,k_2^2)^{\frac{1}{2}}
(h_2h_3)^{\theta_3+\varepsilon}}
{k_2^{\frac{5}{2}-\varepsilon}}
(1+|\frac{3s}{2}|+|w+\frac{s}{2}|)
^{\frac{3}{4}+2\varepsilon}\\
&\hskip 30pt \times(3+|w+\frac{s}{2}|+|\frac{s}{2}|
+\sum_{i=1}^3|\gamma_i|)^{\frac{3}{4}
+\varepsilon}
|\check{\Phi}(\frac{s}{2}+w)|
|\prod_{i=1}^3\Gamma(\frac{s+\frac{1}{2}
+\gamma_i}{2})||\frac{e^{s^2}}{s}| |\dd w \dd s|
\\
&\ll_\Phi \frac{l^{\frac{3}{4}+3\varepsilon}
Y^{\frac{3}{2}+2\varepsilon}}
{X^{\frac{1}{8}-4\varepsilon}},
\end{split}
\end{equation}
because $(l,k_2^2)^{\frac{1}{2}}
h_2^{\theta_3+\varepsilon}\leq (k_2,l^{\infty})^{1+\frac{5}{14}+\varepsilon}$.

\subsection{The contribution of the remainder terms $S_R(\chi_{8d}(l)B(d);\Phi)$}
Observe that $R_Y(d)$ equals $0$ unless $d=d_1 d_2^2$ where $d_1$ is square-free and $d_2>Y$.

Hence
\begin{equation*}
S_R(\chi_{8d}(l)B(d);\Phi)= X^{-1}\sum_{\substack{Y<d_2\leq \sqrt{2X}\\ (d_2,2l)=1}}\mu(d_2)
\sum_{(d_1,l)=1}\nolimits^\flat \chi_{8d_1}(l)B(d_1d_2^2)\Phi(\frac{d_1}{X/d_2^2}),
\end{equation*}
and
\begin{equation*}
B(d_1d_2^2)=\frac{1}{2\pi i}\int_{(u)}\prod_{j=1}^3\frac{
\Gamma(\frac{s+\frac{1}{2}-\gamma_j}{2})}
{\Gamma(\frac{\frac{1}{2}-\gamma_j}{2})}
(\frac{8d_1d_2^2}{\pi})^{\frac{3}{2}s}
\sum_{n=1}^\infty \chi_{8d_1d_2^2}(n)\frac{A(n,1)}
{n^{s+\frac{1}{2}}}e^{s^2}\frac{\dd s}{s}.
\end{equation*}
Plainly
\begin{equation*}
\sum_{n=1}^\infty \chi_{8d_1d_2^2}(n)\frac{A(n,1)}{n^{s+\frac{1}{2}}}
=L(\frac{1}{2}+s,\phi\otimes\chi_{8d_1})
I_{d_1}(s,d_2),
\end{equation*}
where
\begin{equation*}
I_{d_1}(s,d_2)=\prod_{p\mid d_2}(1-\frac{A(p,1)\chi_{8d_1}(p)}{p^{s+\frac{1}{2}}}
+\frac{A(1,p)\chi_{8d_1}(p)^2}{p^{2s+1}}
-\frac{\chi_{8d_1}(p)}{p^{3s+\frac{3}{2}}}).
\end{equation*}
Hence we may move the line of integration to the line $\Re (s)=\frac{1}{\log X}$.
This gives
\begin{equation*}
B(d_1d_2^2)=\frac{1}{2\pi i}\int_{(\frac{1}{\log X})}\prod_{j=1}^3\frac{
\Gamma(\frac{s+\frac{1}{2}-\gamma_j}{2})}
{\Gamma(\frac{\frac{1}{2}-\gamma_j}{2})}
(\frac{8d_1d_2^2}{\pi})^{\frac{3}{2}s}
L(\frac{1}{2}+s,\phi\otimes\chi_{8d_1})
I_{d_1}(s,d_2)
e^{s^2}\frac{\dd s}{s}.
\end{equation*}
Plainly
\begin{equation*}
|I_{d_1}(s,d_2)|\ll t^{\varepsilon}\prod_{p|d_2}
(1+p^{\theta_3-\frac{1}{2}- \Re s}).
\end{equation*}
From $e^{s^2}$ is a Schwartz function,
we truncate the integral from $\frac{1}{\log X}-X^\varepsilon i$ to $\frac{1}{\log X}+X^\varepsilon i$,
Thus
\begin{equation*}
B(d_1d_2^2)=\frac{1}{2\pi i}\int_{\frac{1}{\log X}-X^\varepsilon i}
^{\frac{1}{\log X}+X^\varepsilon i} \prod_{j=1}^3\frac{
\Gamma(\frac{s+\frac{1}{2}-\gamma_j}{2})}
{\Gamma(\frac{\frac{1}{2}-\gamma_j}{2})}
(\frac{8d_1d_2^2}{\pi})^{\frac{3}{2}s}
L(\frac{1}{2}+s,\phi\otimes\chi_{8d_1})
I_{d_1}(s,d_2)
e^{s^2}\frac{\dd s}{s}+O(\frac{1}{X^{2022}}),
\end{equation*}
and
\begin{equation}\label{eqn:finalSR<<1}
\begin{split}
S_R(\chi_{8d}(l)B(d);\Phi)&= \frac{1}{2\pi i X}\int_{\frac{1}{\log X}-X^\varepsilon i}
^{\frac{1}{\log X}+X^\varepsilon i}
\prod_{j=1}^3
\frac{\Gamma(\frac{s+\frac{1}{2}-\gamma_j}{2})}
{\Gamma(\frac{\frac{1}{2}-\gamma_j}{2})}\\
     &\hskip -2cm  \times
\sum_{\substack{Y<d_2\leq \sqrt{2X}\\ (d_2,2l)=1}}\mu(d_2)d_2^{3s} \sum_{d_3\mid d_2}\frac{\mu(d_3)A(d_3,1)}
     {d_3^{s+\frac{1}{2}}}
     \sum_{d_4\mid\frac{d_2}{d_3}}\frac{A(1,d_4)}
     {d_4^{2s+1}}
     \sum_{d_5\mid\frac{d_2}{d_3d_4}}
     \frac{1}{d_5^{3s+\frac{3}{2}}}
     \\
     &\hskip -2cm \times \sum_{d_1}\nolimits^\flat \chi_{8d_1}(ld_3d_5d_4^2)
     (\frac{8d_1}{\pi})^{\frac{3}{2}s}
     L(\frac{1}{2}+s,\phi\otimes\chi_{8d_1})
     \Phi(\frac{d_1}{X/d_2^2})e^{s^2}\frac{\dd s}{s}
+O(\frac{1}{X^{2022}}).
\end{split}
\end{equation}
By using the approximate functional equations of $L(\frac{1}{2}+s_0,\phi\otimes\chi_{8d})$ and the large sieve estimate for quadratic characters in \cite{HB95},
we have
\begin{equation}\label{eqn:5/4}
|\sum_{d}\nolimits^\flat \chi_{8d}(l_0)
(\frac{8d}{\pi})^{\frac{3}{2}s_0} L(\frac{1}{2}+s_0,\phi\otimes\chi_{8d})
\Phi(\frac{d}{Z})|\ll
Z^{\frac{5}{4}+\varepsilon}|1+s_0|^B
l_0^{\varepsilon}
\end{equation}
for some constant $B>0$.
Then we truncate the sum of $d_2$ in \eqref{eqn:finalSR<<1} from $Y$ to $X^{\frac{1}{3}}$, getting
\begin{equation}\label{eqn:finalSR<<}
\begin{split}
S_R(\chi_{8d}(l)B(d);\Phi)&= \frac{1}{2\pi i X}\int_{\frac{1}{\log X}-X^\varepsilon i}
^{\frac{1}{\log X}+X^\varepsilon i}
\prod_{j=1}^3
\frac{\Gamma(\frac{s+\frac{1}{2}-\gamma_j}{2})}
{\Gamma(\frac{\frac{1}{2}-\gamma_j}{2})}\\
     &\hskip -2cm  \times
\sum_{\substack{Y<d_2\leq X^{\frac{1}{3}}\\ (d_2,2l)=1}}\mu(d_2)d_2^{3s} \sum_{d_3\mid d_2}\frac{\mu(d_3)A(d_3,1)}
     {d_3^{s+\frac{1}{2}}}
     \sum_{d_4\mid\frac{d_2}{d_3}}\frac{A(1,d_4)}
     {d_4^{2s+1}}
     \sum_{d_5\mid\frac{d_2}{d_3d_4}}
     \frac{1}{d_5^{3s+\frac{3}{2}}}
     \\
     &\hskip -2cm \times \sum_{d_1}\nolimits^\flat \chi_{8d_1}(ld_3d_5d_4^2)
     (\frac{8d_1}{\pi})^{\frac{3}{2}s}
     L(\frac{1}{2}+s,\phi\otimes\chi_{8d_1})
     \Phi(\frac{d_1}{X/d_2^2})e^{s^2}\frac{\dd s}{s}
+O( X^{-\frac{1}{4}+\varepsilon}).
\end{split}
\end{equation}
But we can not get strong enough estimate of $S_R$ by simply applying \eqref{eqn:5/4} to all $d_2$.
We should use a recursive argument as in Young \cite{Young} to prove the following lemma.

\begin{lemma} \label{lemma:SR}
  Assume $1\leq Y\leq X^{1/4}$. Then
   we have
\begin{equation}\label{eqn:SR<<}
S_R(\chi_{8d}(l)B(d);\Phi)
\ll l^{\frac{3}{4}+\varepsilon}
\frac{X^{\varepsilon}}{Y}.
\end{equation}
\end{lemma}

\begin{proof}


We would show that if for all $Z\leq 2X$ and $\frac{1}{\log X}\leq \Re(s_0)\leq \frac{1000}{\log X}$, $-X^{\varepsilon}\leq \Im(s_0)\leq X^{\varepsilon}$, we have
\begin{equation}\label{iterate1}
|\sum_{d}\nolimits^\flat \chi_{8d}(l_0)
(\frac{8d}{\pi})^{\frac{3}{2}s_0} L(\frac{1}{2}+s_0,\phi\otimes\chi_{8d})
\Phi(\frac{d}{Z})|\ll
l_0^{\varepsilon}Z^{1+\delta}|1+s_0|^B
+l_0^{\frac{3}{4}+\varepsilon}
Z^{\frac{19}{20}+\varepsilon}|1+s_0|^B
\end{equation}
with some constants $\delta\geq 0$ and $B>0$,
then we have
\begin{equation}\label{iterate2}
 \begin{split}
|\sum_{d}\nolimits^\flat \chi_{8d}(l_0)
(\frac{8d}{\pi})^{\frac{3}{2}s_0} & L(\frac{1}{2}+s_0-\frac{1}{\log X}+iv,\phi\otimes\chi_{8d})
\Phi(\frac{d}{Z})| \\
& \ll
(l_0^{\varepsilon}Z^{1+\delta-\frac{1}{20}}
+l_0^{\varepsilon}Z^{1+\varepsilon}
+l_0^{\frac{3}{4}+\varepsilon}
Z^{\frac{19}{20}+\varepsilon})
|1+s_0|^{B'} (1+|v|)^{B''}
\end{split}
\end{equation}
with some constants $B'\geq B$ and $B''>0$.

Now we could iterate from $\delta=\frac{1}{4}+\varepsilon$ as in \eqref{eqn:5/4}, set $\Re (s_0)=\frac{6}{\log X}$,
and iterate five times, then we have
\begin{multline*}
\sum_{d_1}\nolimits^\flat \chi_{8d_1}(ld_3d_5d_4^2)
     (\frac{8d_1}{\pi})^{\frac{3}{2}s} L(\frac{1}{2}+s,\phi\otimes\chi_{8d_1})
     \Phi(\frac{d_1}{X/d_2^2}) \\
 \ll ((ld_3d_5d_4^2)^{\frac{3}{4}+\varepsilon}
     X^{\frac{19}{20}+\varepsilon}
     d_2^{-\frac{19}{10}+\varepsilon}
     +l^{\varepsilon}X^{1+\varepsilon}
     d_2^{-2+\varepsilon})|1+s|^A
\end{multline*}
with some constant $A$.
Now by \eqref{eqn:finalSR<<} we prove \eqref{eqn:SR<<}.


Now we still need to prove \eqref{iterate2} by using \eqref{iterate1}.
As in \eqref{eqn:sum2S}, we need to calculate
\begin{equation*}
S_1=\sum_{n=1}^{\infty}\frac{A(n,1)}{n^{\frac{1}{2}+s_0}}
\sum_{d}\nolimits^\flat \chi_{8d}(l_0n)
(\frac{8d}{\pi})^{\frac{3}{2}s_0} \frac{1}{2\pi i}\int_{(u)}\prod_{j=1}^{3}
\frac{\Gamma(\frac{s+\frac{1}{2}+s_0-\gamma_j}{2})}
{\Gamma(\frac{\frac{1}{2}+s_0-\gamma_j}{2})}
(n(\frac{\pi}{8d})^{\frac{3}{2}})^{-s}
e^{s^2}\frac{\dd s}{s}
\Phi(\frac{d}{Z})
\end{equation*}
and
\begin{equation*}
S_2=\sum_{n=1}^{\infty}\frac{A(1,n)}{n^{\frac{1}{2}-s_0}}
\sum_{d}\nolimits^\flat \chi_{8d}(l_0n)
(\frac{8d}{\pi})^{-\frac{3}{2}s_0} \frac{1}{2\pi i}\int_{(u)}\prod_{j=1}^{3}
\frac{\Gamma(\frac{s+\frac{1}{2}-s_0+\gamma_j}{2})}
{\Gamma(\frac{\frac{1}{2}+s_0-\gamma_j}{2})}
(n(\frac{\pi}{8d})^{\frac{3}{2}})^{-s}
e^{s^2}\frac{\dd s}{s}
\Phi(\frac{d}{Z}).
\end{equation*}
We also use $\mu^2(d)=M_{Y_Z}(d)+R_{Y_Z}(d)$, and we have
\begin{equation*}
  S_1=Z S_{1,M}+Z S_{1,R}
\end{equation*}
where
\begin{equation*}
  S_{1,M}=
  \sum_{n=1}^{\infty}
  \frac{A(n,1)}{n^{\frac{1}{2}+s}}
  S_{M,Z,Y_{Z}}(\chi_{8d}(l_0n);\Phi_{1,n,s}),
\end{equation*}
\begin{equation*}
  S_{1,R}=
  \sum_{n=1}^{\infty}
  \frac{A(n,1)}{n^{\frac{1}{2}+s}}
  S_{R,Z,Y_{Z}}(\chi_{8d}(l_0n);\Phi_{1,n,s}),
\end{equation*}
with $S_{M,Z,Y_{Z}}$ and $S_{R,Z,Y_{Z}}$ in \S2, and
\begin{equation*}
\Phi_{1,y,s}(t)=\frac{\Phi(t)}{2\pi i}(\frac{8Zt}{\pi})^{\frac{3s_0}{2}}
\int_{(u)}\prod_{j=1}^{3}
\frac{\Gamma(\frac{s+\frac{1}{2}-s_0+\gamma_j}{2})}
{\Gamma(\frac{\frac{1}{2}+s_0-\gamma_j}{2})}
(y(\frac{\pi}{8Zt})^{\frac{3}{2}})^{-s}
e^{s^2}\frac{\dd s}{s}.
\end{equation*}
Then by Poisson summation formula we have $S_{1,M}=P_1+R_1$ according to $k=0$ and $k\neq 0$,
then using the same method as \S \ref{mainterm},
from
\begin{equation*}
  \frac{1}{\Gamma(\frac{\frac{1}{2}+s-\gamma_j}
  {2})}\ll \Im(s)^{\frac{1}{2}},
\end{equation*}
and
\begin{equation*}
  G_\phi(1+s_0;l_0)\ll l_0^{\frac{1}{2}+\varepsilon},
\end{equation*}
we have $P_1\ll l_0^{\varepsilon}Z^{\varepsilon}
+l_0^{\frac{1}{2}}
Z^{-\frac{3}{14}+\varepsilon}
Y_{Z}^{\frac{3}{14}+\varepsilon}$.
Also by using the same method as \S \ref{remainder1}, we have $R_1\ll l_0^{\frac{3}{4}+\varepsilon}
Y_{Z}^{\frac{3}{2}+\varepsilon}
Z^{-\frac{1}{8}+\varepsilon}$.
By estimating terms in $S_{1,R}$ with $Y_{Z}<d_2\leq Z^{\frac{1}{3}}$ like our previous discussion about \eqref{eqn:finalSR<<} with \eqref{iterate1},  estimating terms in $S_{1,R}$ with $d_2> Z^{\frac{1}{3}}$ with \eqref{eqn:5/4}, and taking $Y_{Z}=Z^{1/20}$, we have
\begin{equation*}
\begin{split}
S_{1,R}(\chi_{8d}(l)B(d);\Phi)&\ll \frac{1}{Z}\int_{\frac{1}{\log X}-Z^\varepsilon i}
^{\frac{1}{\log X}+Z^\varepsilon i}
\prod_{j=1}^3
|\frac{\Gamma(\frac{s+s_0+\frac{1}{2}-\gamma_j}{2})}
{\Gamma(\frac{\frac{1}{2}-\gamma_j}{2})}|
\sum_{\substack{Y_{Z}<d_2\leq Z^{\frac{1}{3}}\\ (d_2,2l_0)=1}}|\mu(d_2)|\sum_{d_3\mid d_2}
     d_3^{-\frac{1}{2}+\theta_3}\\
     &\hskip 10pt  \times
     \sum_{d_4\mid\frac{d_2}{d_3}}
     d_4^{-1+\theta_3}
     \sum_{d_5\mid\frac{d_2}{d_3d_4}}
     d_5^{-\frac{3}{2}}
     \Big((l_0d_3d_5d_4^2)^{\varepsilon}
     (\frac{Z}{d_2^2})^{1+\delta}
+(l_0d_3d_5d_4^2)^{\frac{3}{4}+\varepsilon}
(\frac{Z}{d_2^2})^{\frac{19}{20}+\varepsilon}
\Big)     \\
     &\hskip 60pt \times
|1+s_0+s|^B
     |e^{s^2}|
     \frac{\dd s}{|s|}
+O( Z^{-\frac{1}{4}+\varepsilon})\\
&\ll l_0^{\frac{3}{4}+\varepsilon}
Z^{-\frac{1}{20}+\varepsilon}
+l_0^{\varepsilon}
Z^{\delta-\frac{1}{20}}.
\end{split}
\end{equation*}
Hence we get
\[
  S_1 \ll l_0^{\varepsilon}Z^{1+\varepsilon}
  + l_0^{\frac{3}{4}+\varepsilon} Z^{\frac{19}{20}+\varepsilon}
  +l_0^{\varepsilon} Z^{1+\delta-\frac{1}{20}}.
\]
Similarly for $S_2$. This proves \eqref{iterate2},
and  completes the proof of the lemma.
\end{proof}

Now by \eqref{eqn:sum2S}, \eqref{eqn:P=1}, \eqref{eqn:P=2}, \eqref{eqn:R<<}, and Lemma \ref{lemma:SR}, if we take $Y=X^{1/20}$,   we have when $\phi$ is not self-dual
\begin{multline}
\sum_{2\nmid d}\nolimits^\flat L(\frac{1}{2},\phi\otimes \chi_{8d})\chi_{8d}(l)\Phi(\frac{d}{X})
=\frac{2\check{\Phi}(0)X}{3\zeta(2)\sqrt{l_1}}
\prod_{p\mid l}\frac{p}{p+1}
 \bigg(G_\phi(l)L^{\{2\}}(1,\sym^2 \phi) \\ +\prod_{i=1}^3\frac{\Gamma(\frac{
 \frac{1}{2}+\gamma_i}{2})}
 {\Gamma(\frac{\frac{1}{2}-\gamma_i}{2})}
 \bar{G}_\phi(l)L^{\{2\}}
 (1,\sym^2 \widetilde{\phi})\bigg)
+O_\Phi ( l^{\frac{3}{4}+\varepsilon}X^{\frac{19}{20}+\varepsilon}),
\end{multline}
and when $\phi$ is self-dual
\begin{multline}
\sum_{2\nmid d}\nolimits^\flat L(\frac{1}{2},\phi\otimes \chi_{8d})\chi_{8d}(l)\Phi(\frac{d}{X})
= \frac{\lim_{s\to 1} (s-1)L^{\{2\}}(s,\sym^2 \phi)\check{\Phi}(0)}{\zeta(2)\sqrt{l_1}}
\prod_{p\mid l}\frac{p}{p+1}X\\
\times \Big(G_\phi(l)
\log \frac{X}{l_1 ^{\frac{2}{3}}}+C_\phi(l)\Big)
+O_\Phi ( l^{\frac{3}{4}+\varepsilon}
X^{\frac{19}{20}+\varepsilon}).
\end{multline}
This completes the proof of Theorem \ref{thm2}.

\section{Proof of Theorem \ref{thm1}}\label{sec:det_form}

\begin{lemma}\label{lm5.1}
  Let $M$ be a fixed integer coprime with $3,5,7,11$,
  $G_\phi(l)$ be defined as in Theorem \ref{thm2}.
  Let  $Q_\phi$ be a positive integer such that $(Q_{\phi},2M)=1$ and $G_\phi(Q_{\phi})\neq 0$. For normalized $\phi$ and $\phi'$, if we have $G_{\phi}(Q_{\phi}N)=a \cdot G_{\phi'}(Q_{\phi}N)$ for all $N$ coprime to $2Q_{\phi}M$, with some nonzero constant $a$, then $\phi=\phi'$.
\end{lemma}

\begin{proof}
By Theorem \ref{thm2} we know $G_{\phi}(c_\phi H^2)\neq 0$ for all $(H,2c_\phi)=1$.
So there exists a positive integer  $Q_\phi$ such that $(Q_{\phi},2M)=1$ and $G_\phi(Q_{\phi})\neq 0$.
Let $A(m,n)$ be the coefficients of $\phi$ and $A'(m,n)$ be the coefficients of $\phi'$.
From $G_{\phi}(Q_\phi)=a \cdot G_{\phi'}(Q_\phi)$ with $a\neq 0$,
we have $G_{\phi'}(Q_\phi)\neq 0$.
For any odd prime $p\geq 13$ satisfying $(p,2Q_\phi M)=1$, by comparing both sides of  $G_{\phi}(Q_\phi N)=a \cdot G_{\phi'}(Q_\phi N)$ with $N=p$ and $p^2$,
and the fact $G_{\phi}(Q_\phi p^2),G_{\phi'}(Q_\phi p^2)\neq 0$,
we have
\begin{equation*}
G_{\phi}(Q_\phi p)/G_{\phi}(Q_\phi p^2)=
G_{\phi'}(Q_\phi p)/G_{\phi'}(Q_\phi p^2).
\end{equation*}
Hence we have
\begin{equation*}
G_{\phi,p}(Q_\phi p)/G_{\phi,p}(Q_\phi p^2)=
G_{\phi',p}(Q_\phi p)/G_{\phi',p}(Q_\phi p^2).
\end{equation*}
From the definition \eqref{eqdefineta},
we have
\begin{equation*}
\frac{1+pA(p,1)}{p+A(1,p)}=
\frac{1+pA'(p,1)}{p+A'(1,p)},
\end{equation*}
and hence
\begin{equation}\label{eqn:Aphi==}
p^2A(p,1)+pA(p,1)\overline{A'(p,1)}
+\overline{A'(p,1)}=
p^2A'(p,1)+pA'(p,1)\overline{A(p,1)}
+\overline{A(p,1)}.
\end{equation}
By comparing the real parts of both sides, we get
\begin{multline*}
p^2\Re(A(p,1))+p\Re(A(p,1))\Re(A'(p,1))
+p\Im(A(p,1))\Im(A'(p,1))+\Re(A'(p,1))
\\=p^2\Re(A'(p,1))+p\Re(A(p,1))\Re(A'(p,1))
+p\Im(A(p,1))\Im(A'(p,1))+\Re(A(p,1)).
\end{multline*}
So we get $(p^2-1)\Re(A(p,1))
=(p^2-1)\Re(A'(p,1))$, from which  we obtain
\begin{equation}\label{realpart}
\Re(A(p,1))=\Re(A'(p,1)).
\end{equation}
Then by comparing the imaginary parts of both sides in \eqref{eqn:Aphi==}, we get
\begin{multline*}
p^2\Im(A(p,1))-p\Re(A(p,1))\Im(A'(p,1))
+p\Re(A'(p,1))\Im(A(p,1))-\Im(A'(p,1))
\\=p^2\Im(A'(p,1))+p\Re(A(p,1))\Im(A'(p,1))
-p\Re(A'(p,1))\Im(A(p,1))
-\Im(A(p,1)).
\end{multline*}
Together with \eqref{realpart}, we have
\begin{equation}
  (p^2+2p\Re(A(p,1))+1)
  (\Im(A(p,1))-\Im(A'(p,1)))=0.
\end{equation}
By the bounds toward to the Ramanujan conjecture, we get   $|\Re(A(p,1))| \leq |A(p,1)| \leq 3 p^{\frac{5}{14}} < \frac{p^2+1}{2p}$, and hence $p^2+2p\Re(A(p,1))+1 > 0$, provided $p\geq 17$.
Thus for $(p,2Q_\phi M)=1$ and $p\geq 17$, we have
\begin{equation}
\Im(A(p,1))=\Im(A'(p,1))
\end{equation}
Thus we know that $A(p,1)=A'(p,1)$ for all odd primes $p\geq 17$ which are coprime to $2Q_\phi M$.
By the strong multiplicity one theorem (see e.g. \cite[\S12.6]{GB06}), we prove $\phi=\phi'$.
\end{proof}

\begin{proof}[Proof of Theorem \ref{thm1}]

We first claim that the equality of  twisted $L$-values  determines whether the forms are both self-dual or not. By taking $l=c_\phi$ in Theorem \ref{thm2}, we get $G_\phi(l)\neq0$. So we know $S(\phi) := \sum_{2\nmid d}\nolimits^\flat \chi_{8d}(l) L(\frac{1}{2},\phi\otimes \chi_{8d})\Phi(\frac{d}{X})$ is $ \asymp X\log X$ if $\phi$ is self-dual and $\ll X$ if not. Since $L(\frac{1}{2},\phi\otimes \chi_{8d})=\kappa L(\frac{1}{2},\phi'\otimes \chi_{8d})$ for all odd square-free $d$ and $\kappa\neq0$, we get $S(\phi') = \kappa S(\phi)$ is  $ \asymp X\log X$ if $\phi$ is self-dual. By Theorem \ref{thm2} again, we know $S(\phi')\asymp X\log X$ holds only if $\phi'$ is self-dual.
Hence we show that if $\phi$ is self-dual then $\phi'$ is also.
Similarly, we know if $\phi'$ is self-dual then $\phi$ is also. This proves our claim.

Case i).
If $\phi$ and $\phi'$ are both self-dual, then by comparing the main terms in Theorem \ref{thm2} we have
\begin{equation*}
G_{\phi}(l)=\bar{G}_{\phi}(l)=a G_{\phi'}(l)=a \bar{G}_{\phi'}(l)
\end{equation*}
with some nonzero constant $a$ depending on $\phi$ and $\phi'$ but not $l$.
Then in Lemma \ref{lm5.1}
let $Q_\phi=c_\phi$ we have $\phi=\phi'$.

Case ii).
If  $\phi$ and $\phi'$ are both not self-dual,
from Theorem \ref{thm2} we have
\begin{equation} \label{eqn:Gphi=Gphi'}
a_1G_{\phi}(l)+a_2\bar{G}_{\phi}(l)=
b_1G_{\phi'}(l)+b_2 \bar{G}_{\phi'}(l)
\end{equation}
where $a_1=L^{\{2\}}(1,\sym^2 \phi)$,
$a_2=\prod_{i=1}^3\frac{\Gamma(\frac{\frac{1}{2}
 +\gamma_{\phi,i}}{2})}{\Gamma(\frac{\frac{1}{2}
 -\gamma_{\phi,i}}{2})}
L^{\{2\}}(1,\sym^2 \widetilde{\phi})$,
$b_1=L^{\{2\}}(1,\sym^2 \phi')$,
$b_2=\prod_{i=1}^3\frac{\Gamma(\frac{\frac{1}{2}
 +\gamma_{\phi',i}}{2})}{\Gamma(\frac{\frac{1}{2}
 -\gamma_{\phi',i}}{2})}
L^{\{2\}}(1,\sym^2 \widetilde{\phi'})$ are all non-zero.

By strong multiplicity one theorem and $\phi \neq \widetilde{\phi}$, we know there are infinity primes $p_0\geq 13$ such that $A(p_0,1)\neq A(1,p_0)$.
For such $p_0$, we have $\frac{1+p_0A(p_0,1)}{p_0+A(1,p_0)}\notin
\mathbb{R}$, and $G_{\phi}(c_\phi p_0)\neq 0$.
Fix such a prime $p_0$ with $(p_0,M)=1$.
Then there exists a $d_\phi\in \{c_\phi p_0,c_\phi p_0^2\}$ such that $a_1G_{\phi}(d_\phi)+a_2\bar{G}_{\phi}(d_\phi)\neq 0$.
In particular, we have $G_{\phi}(d_\phi)\neq0$.
Indeed,
if not, then we have
$a_1G_{\phi}(c_\phi p_0)+a_2\bar{G}_{\phi}(c_\phi p_0)
= a_1G_{\phi}(c_\phi p_0^2)+a_2\bar{G}_{\phi}(c_\phi p_0^2)= 0$.
Together with $G_{\phi}(c_\phi p_0)=\frac{1+p_0A(p_0,1)}{p_0+A(1,p_0)}G_{\phi}(c_\phi p_0^2)$ and $\frac{1+p_0A(p_0,1)}{p_0+A(1,p_0)}\notin
\mathbb{R}$, we get $G_{\phi}(c_\phi p_0^2)=0$, which
contradicts with Theorem \ref{thm2}.
Note that by \eqref{eqn:Gphi=Gphi'} we have $b_1G_{\phi'}(d_\phi)+b_2 \bar{G}_{\phi'}(d_\phi)\neq 0$.
So $G_{\phi'}(d_\phi)\neq 0$.

If $\phi=\phi'$, then we finish the proof.

If $\phi\neq \phi'$,  then by the strong multiplicity one theorem (see e.g. \cite[\S12.6]{GB06}), there are infinity primes $p\geq 17$ satisfying $A(p,1)\neq A'(p,1)$. By the same argument as in Lemma \ref{lm5.1} we have  $\frac{1+pA(p,1)}{p+A(1,p)}\neq
\frac{1+pA'(p,1)}{p+A'(1,p)}$ for such $p$'s. Hence for those $p$, we have
\begin{equation}\label{eqn:notequal}
G_{\phi,p}(pd_\phi)/G_{\phi,p}(p^2d_\phi)
\neq G_{\phi',p}(pd_\phi)/G_{\phi',p}(p^2d_\phi).
\end{equation}
We fix one such $p\geq 17$ with $(p,d_\phi M)=1$. By \eqref{eqn:Gphi=Gphi'} with $l=pd_\phi N$ and $l=p^2 d_\phi N$, we obtain
\begin{equation*}
a_1G_{\phi}(p d_{\phi} N)
+a_2\bar{G}_{\phi}(pd_{\phi}N)
=b_1G_{\phi'}(pd_{\phi}N)
+b_2\bar{G}_{\phi'}(pd_{\phi}N),
\end{equation*}
and
\begin{equation*}
a_1G_{\phi}(p^2d_{\phi}N)
+a_2\bar{G}_{\phi}(p^2d_{\phi}N)
=b_1G_{\phi'}(p^2d_{\phi}N)
+b_2\bar{G}_{\phi'}(p^2d_{\phi}N)
\end{equation*}
for any integer $N$ satisfying $(N,2d_{\phi} pM)=1$.
Thus by a linear combination of the above two identities to eliminate $a_1$,
and $G_{\phi}(pd_{\phi}N)=
\frac{G_{\phi,p}(pd_{\phi})}
{G_{\phi,p}(p^2d_{\phi})}
G_{\phi}(p^2d_{\phi}N)$,
we get
\begin{equation}\label{eqn:p2=}
\begin{split}
&(\frac{G_{\phi,p}(pd_{\phi})}
{G_{\phi,p}(p^2d_{\phi})}
-\frac{\bar{G}_{\phi,p}(pd_{\phi})}
{\bar{G}_{\phi,p}(p^2d_{\phi})})
a_2\bar{G}_{\phi}(p^2d_{\phi}N)\\
&=(\frac{G_{\phi,p}(pd_{\phi})}
{G_{\phi,p}(p^2d_{\phi})}
-\frac{G_{\phi',p}(pd_{\phi})}
{G_{\phi',p}(p^2d_{\phi})})
b_1G_{\phi'}(p^2d_{\phi}N)
+(\frac{G_{\phi,p}(pd_{\phi})}
{G_{\phi,p}(p^2d_{\phi})}
-\frac{\bar{G}_{\phi',p}(pd_{\phi})}
{\bar{G}_{\phi',p}(p^2d_{\phi})})
b_2\bar{G}_{\phi'}(p^2d_{\phi}N).
\end{split}
\end{equation}

By \eqref{eqn:notequal}, we  get
\begin{equation}\label{eqn:nonzero}
  \frac{G_{\phi,p}(pd_{\phi})}
  {G_{\phi,p}(p^2d_{\phi})}
-\frac{\bar{G}_{\phi,p}(pd_{\phi})}
{\bar{G}_{\phi,p}(p^2d_{\phi})}\neq 0.
\end{equation}
Indeed, if \eqref{eqn:nonzero} is not true, then by \eqref{eqn:p2=} we have
\begin{equation*}
G_{\phi'}(p^2d_{\phi}N)=a \cdot \bar{G}_{\phi'}(p^2d_{\phi}N)
\end{equation*}
for all $(N,2d_{\phi} pM)=1$ with some constant $a$.
Let $N$ be a square-full number then we must have $a\neq 0$.
Then by Lemma \ref{lm5.1} with $Q_{\phi'}=p^2 d_{\phi}$ we have $\phi'=\widetilde{\phi'}$, which contradicts to that $\phi'$ is not self-dual.

By \eqref{eqn:p2=} and \eqref{eqn:nonzero}, we can rewrite as
\begin{equation}\label{eqn:c2}
\bar{G}_{\phi}(p^2d_{\phi}N)=
c_1G_{\phi'}(p^2d_{\phi}N)
+c_2\bar{G}_{\phi'}(p^2d_{\phi}N)
\end{equation}
for all $(N,2d_{\phi}pM)=1$ with some constants $c_1,c_2$ independent of $N$, and we know $c_1\neq 0$.
Note that as in \eqref{eqn:notequal}, there are infinity many primes $q$ with $(q,pM)=1$ such that
\begin{equation*}
G_{\phi,q}(qp^2d_{\phi})/
G_{\phi,q}(q^2p^2d_{\phi})\neq
G_{\phi',q}(qp^2d_{\phi})/
G_{\phi',q}(q^2p^2d_{\phi}).
\end{equation*}
We fix one such $q$. By similar arguments as above we can eliminate $c_2$ in \eqref{eqn:c2}, getting
\begin{multline*}
(\frac{G_{\phi',q}(qp^2d_{\phi})}
{G_{\phi',q}(q^2p^2d_{\phi})}
-\frac{G_{\phi,q}(qp^2d_{\phi})}
{G_{\phi,q}(q^2p^2d_{\phi})})
G_{\phi}(q^2p^2d_{\phi}N)
\\
=
(\frac{G_{\phi',q}(qp^2d_{\phi})}
{G_{\phi',q}(q^2p^2d_{\phi})}
-\frac{\bar{G}_{\phi',q}(qp^2d_{\phi})}
{\bar{G}_{\phi',q}(q^2p^2d_{\phi})})
\bar{c}_1\bar{G}_{\phi'}(q^2p^2d_{\phi}N),
\end{multline*}
for all $(N,2pqd_{\phi}M)=1$.
Let $N$ be a square-full number then we know $(\frac{G_{\phi',p}(qp^2d_{\phi})}
{G_{\phi',p}(q^2p^2d_{\phi})}
-\frac{\bar{G}_{\phi',p}(qp^2d_{\phi})}
{\bar{G}_{\phi',p}(q^2p^2d_{\phi})})\bar{c}_1\neq 0$.
Thus we have $G_{\phi}(q^2p^2d_{\phi}N)=a \cdot \bar{G}_{\phi'}(q^2p^2d_{\phi}N)$ for all $(N,2pqd_{\phi}M)=1$ with some nonzero constant $a$ depends on $\phi,\phi',p,q,M,d_\phi$.
Then by Lemma \ref{lm5.1} with $Q_\phi=q^2p^2d_\phi$, we have $\phi=\widetilde{\phi'}$.

In conclusion, if $\phi$ and $\phi'$ are both not self-dual, then we have
$\phi= \phi'\textrm{ or }\widetilde{\phi'}.$
This completes the proof of Theorem \ref{thm1}.
\end{proof}

\section{Proof of Theorem \ref{thm3}}\label{sec:FC}
\begin{proof}[Proof of Theorem \ref{thm3}]
By the Rankin--Selberg theory, we have  (see \cite{Mo02})
\begin{equation*}
\sum_{m^2n\leq X}|A(m,n)|^2\ll X^{1+\varepsilon},
\end{equation*}
and hence
\begin{equation} \label{eqn:RS11}
\sum_{n\leq X}|A(n,1)|\ll X^{1+\varepsilon},
\end{equation}
Recall that
\begin{equation*}
L(s,\phi\otimes  \phi)=\mathop{\sum\sum\sum}_{k,m,n \geq1}\frac{A(m,n)^2}{(k^3m^2n)^s}, \quad  \Re(s)>1.
\end{equation*}
Denote
\begin{equation*}
\lambda_{\phi\otimes \phi}(q)=\sum_{k^3m^2n=q} A(m,n)^2.
\end{equation*}
Then we have
\begin{equation}\label{eqn:RS22}
\sum_{q\leq X} |\lambda_{\phi\otimes \phi}(q)|  \leq \sum_{q\leq X} \lambda_{\phi\otimes \widetilde \phi}(q)  \ll X^{1+\varepsilon}.
\end{equation}

Denote $L_p(s,\sym^2 \phi)=\sum_{h=0}^{\infty}\frac{B(p^h,1)}{p^{sh}}$ and $L_p(s,\sym^2 \phi)/L_p(s,\widetilde \phi)=\sum_{h=0}^{\infty}\frac{C(p^h,1)}{p^{sh}}$.
For prime $p$, let $A(p^{h},1)=B(p^{h},1)=C(p^{h},1)=0$ if $h< 0$.
From
\[
  L_p(s,\sym^2 \phi) = \frac{L_p(s,\phi\otimes \phi)}{ L_p(s,\widetilde{\phi})}= \sum_{h\geq0} \frac{\lambda_{\phi\otimes \phi}(p^h)}{p^{hs}} (1-\frac{A(1,p)}{p^{s}}+\frac{A(p,1)}{p^{2s}}-\frac{1}{p^{3s}})
\]
we have
\begin{equation}\label{eqeqB}
B(p^{h},1)=\lambda_{\phi\otimes \phi}(p^{h})-A(1,p)\lambda_{\phi\otimes \phi}(p^{h-1})+A(p,1)\lambda_{\phi\otimes \phi}(p^{h-2})-\lambda_{\phi\otimes \phi}(p^{h-3}).
\end{equation}
Similarly we have
\begin{equation}\label{eqeqC}
  C(p^{h},1)=
  B(p^{h},1)-A(1,p)B(p^{h-1},1)
  +A(p,1)B(p^{h-2},1)-B(p^{h-3},1),
\end{equation}
and  by \eqref{eqn:sum-2h} we have
\begin{equation}\label{eqeqD}
  A(p^{2h},1)=C(p^{h},1)+A(1,p)C(p^{h-1},1).
\end{equation}

Define $B(n,1)$ by multiplicativity. We  first  prove
\begin{equation}\label{eqn:sumB<<}
  \sum_{n\leq X}|B(n,1)|\ll X^{1+\varepsilon}
\end{equation}
based on \eqref{eqn:RS22} and \eqref{eqeqB}.
By \eqref{eqeqB} we have
\begin{equation*}
\begin{split}
|B(n,1)|\leq \prod_{p^{\alpha}\parallel n}&(
|\lambda_{\phi\otimes \phi}(p^{\alpha})|+|A(p,1)||\lambda_{\phi\otimes \phi}(p^{\alpha-1})|+|A(p,1)||\lambda_{\phi\otimes \phi}(p^{\alpha-2})|  +|\lambda_{\phi\otimes \phi}(p^{\alpha-3})|).
\end{split}
\end{equation*}
Since $|A(n,1)|$ and $|\lambda_{\phi\otimes \phi}(n)|$ are multiplicative, we have
\begin{equation*}
\begin{split}
 \prod_{p^{\alpha} \parallel n}&(
|\lambda_{\phi\otimes \phi}(p^{\alpha})|+|A(p,1)||\lambda_{\phi\otimes \phi}(p^{\alpha-1})|+|A(p,1)||\lambda_{\phi\otimes \phi}(p^{\alpha-2})|  +|\lambda_{\phi\otimes \phi}(p^{\alpha-3})|)
\\ &
= \sum_{\substack{n=n_0n_1n_2n_3 \\ (n_i,n_j)=1, \ i\neq j }} |\lambda_{\phi\otimes \phi}(n_0)|
|A(\textrm{rad }n_1,1)||\lambda_{\phi\otimes \phi}(\frac{n_1}{\textrm{rad }n_1})|
\\ & \hskip 60pt \cdot
|A(\textrm{rad }n_2,1)||\lambda_{\phi\otimes \phi}(\frac{n_2}{(\textrm{rad }n_2)^2})|
|\lambda_{\phi\otimes \phi}(\frac{n_3}{(\textrm{rad }n_3)^3})|.
\end{split}
\end{equation*}
Here $\textrm{rad } n=\prod_{p\mid n}p$ is the radical of $n$.
So  we have
\begin{equation*}
\begin{split}
\sum_{n\leq X}|B(n,1)|& \leq\sum_{n_0\leq X}|\lambda_{\phi\otimes \phi}(n_0)|
\sum_{\substack{n_1\leq X/n_0\\(n_1,n_0)=1}}|A(1,\textrm{rad }n_1)||\lambda_{\phi\otimes \phi}(\frac{n_1}{\textrm{rad }n_1})|\\
&\hskip 15pt \times \sum_{\substack{n_2\leq X/n_0n_1\\(n_2,n_0n_1)=1\\p^2\mid n_2\textrm{ if }p\mid n_2}}
|A(\textrm{rad }n_2,1)||\lambda_{\phi\otimes \phi}(\frac{n_2}{(\textrm{rad }n_2)^2})|
\sum_{\substack{n_3\leq X/n_0n_1n_2\\(n_3,n_0n_1n_2)=1\\p^3\mid n_3\textrm{ if }p\mid n_3}}
|\lambda_{\phi\otimes \phi}(\frac{n_3}{(\textrm{rad }n_3)^3})|
\\
& \leq \sum_{m_1\leq X}|\lambda_{\phi\otimes \phi}(m_1)|
\sum_{m_2\leq X/m_1}|A(1,m_2)|
\sum_{m_3\leq X/m_1m_2}|\lambda_{\phi\otimes \phi}(m_3)|\\
&\hskip 60pt \times
\sum_{m_4\leq X/m_1m_2m_3}|A(m_4,1)|
\sum_{m_5\leq X/m_1m_2m_3m_4^2}|\lambda_{\phi\otimes \phi}(m_5)|\\
& \hskip 60pt \times
\sum_{m_6\leq X/m_1m_2m_3m_4^2m_5}
\sum_{m_7\leq X/m_1m_2m_3m_4^2m_5m_6^3}|\lambda_{\phi\otimes \phi}(m_7)|
\\
& \ll X^{1+\varepsilon}.
\end{split}
\end{equation*}
Here we have used \eqref{eqn:RS11} and \eqref{eqn:RS22}.
This completes the proof of \eqref{eqn:sumB<<}.

Define $C(n,1)$ by multiplicativity.
By \eqref{eqeqC}, \eqref{eqn:sumB<<} and the same process as the proof of \eqref{eqn:sumB<<}, we prove
\begin{equation}\label{eqn:sumC<<}
  \sum_{n\leq X}|C(n,1)|\ll X^{1+\varepsilon}.
\end{equation}
From \eqref{eqeqD} we have
\begin{equation}\label{algo2}
\sum_{n\leq X}|A(n^2,1)|\leq \sum_{n\leq X}\prod_{p^{\alpha}\mid \mid n}(|C(p^{\alpha},1)|+|A(1,p)|
|C(p^{\alpha-1},1)|).
\end{equation}
By \eqref{eqn:sumC<<} and a similar but simpler process as the proof of \eqref{eqn:sumB<<},
we obtain $\sum_{n\leq X}|A(n^2,1)|\ll X^{1+\varepsilon}$,
which completes the proof of Theorem \ref{thm3}.
\end{proof}

\section*{Acknowledgements}
The authors are grateful to the referees for their very helpful comments and suggestions.
They want to thank Prof. Ze\'{e}v Rudnick for his valuable advice. 
They would like to thank Gopal Maiti for pointing out a mistake in the proof of Theorem \ref{thm2} in the previous version of this paper. 


\end{document}